\documentclass[a4paper,11pt]{article} 
\usepackage{amsmath,amssymb,enumerate,color}
\usepackage{amsthm,color}
\usepackage{txfonts}
\usepackage[pagewise]{lineno}
\usepackage{bigints}
\usepackage{comment}
\usepackage{graphicx}

\newcommand{\R}{\mathbb{R}}

\newcommand{\N}{\mathbb{N}}

\newcommand{\norm}[1]{\left\|#1\right\|}
\newcommand{\la}{\langle}
\newcommand{\ra}{\rangle}

\newtheorem{theorem}{Theorem}[section]
\newtheorem{lemma}[theorem]{Lemma}
\newtheorem{proposition}[theorem]{Proposition}

\theoremstyle{remark}
\newtheorem{remark}{Remark}[section]
\theoremstyle{definition}

\newtheorem{definition}{Definition}[section]

\numberwithin{equation}{section}

\makeatletter
\def\@cite#1#2{[{{\bfseries #1}\if@tempswa , #2\fi}]}
\makeatother
\let\oldequation\equation
\let\oldendequation\endequation
\renewenvironment{equation}{\linenomathNonumbers\oldequation}{\oldendequation\endlinenomath}

\begin{document}
\begin{center}
\Large{{\bf
Upper estimates of the lifespan for the  fractional wave equations with time-dependent damping and a power nonlinearity of subcritical and critical Fujita exponent}}
\end{center}

\vspace{5pt}

\begin{center}
Jiayun Lin
\footnote{
Department of Mathematics and Science, School of Sciences, Zhejiang Sci-Tech University, 310018, Hangzhou, China,
E-mail:\ {\tt jylin@zstu.edu.cn}},\
Masahiro Ikeda%
\footnote{
Department of Mathematics, Faculty of Science and Technology, Keio University, 3-14-1 Hiyoshi, Kohoku-ku, Yokohama, 223-8522, Japan/Center for Advanced Intelligence Project, RIKEN, Japan,
E-mail:\ {\tt masahiro.ikeda@keio.jp/masahiro.ikeda@riken.jp}},\

\end{center}

\newenvironment{summary}{\vspace{.5\baselineskip}\begin{list}{}{%
     \setlength{\baselineskip}{0.85\baselineskip}
     \setlength{\topsep}{0pt}
     \setlength{\leftmargin}{12mm}
     \setlength{\rightmargin}{12mm}
     \setlength{\listparindent}{0mm}
     \setlength{\itemindent}{\listparindent}
     \setlength{\parsep}{0pt}
     \item\relax}}{\end{list}\vspace{.5\baselineskip}}
\begin{summary}
{\footnotesize {\bf Abstract.}
In this paper, we investigate the Cauchy problem of the fractional wave equation with time-dependent damping and the source nonlinearity $f(u)\approx |u|^p$:
\begin{equation*}
\begin{cases}
\partial_t^2u(t,x)+(-\Delta)^{\sigma/2} u(t,x)+b(t) \partial_t u(t,x) =f(u(t,x)),\ &(t,x)\ \in [0,T)\times \R^N,\\
u(0,x)=u_0(x),\ \partial_tu(0,x)=u_1(x),\ &x\ \in\ \R^N,
\end{cases}
\end{equation*}
where the damping coefficient behaves as $b(t)\approx (1+t)^{-\beta}$. We first establish local existence and uniqueness of mild and strong solutions in weighted Sobolev spaces by deriving weighted Sobolev estimates for the linear propagators and overcoming the spatial-weight coupling in the nonlinearity through the Bessel potential estimates. Then, in the subcritical and critical cases $1<p\leq p_c:=1+\frac \sigma N$, we derive the upper estimates of the lifespan for the fractional Laplacian with $0<\sigma<2$ and time-dependent damping $\beta \in [-1, 1)$ via the framework of ordinary differential inequalities.
The blow-up results, together with the global existence in the supercritical case obtained in \cite{LI}, show that the critical exponent for the fractional wave equation is $p_c=1+\frac{\sigma}{N}$ for $0<\sigma<2$. Moreover, combining with the lower estimate of lifespan derived in \cite{LI}, we conclude that the estimate in this paper is sharp. Note that our result for the critical case is completely new even in the classical case $b(t)=1$. We also consider the case of $\beta =1$, and obtain the upper estimate of the lifespan.
}
\end{summary}

{\footnotesize{\it Mathematics Subject Classification}\/ (2020): %
35B44; 
35A01; 
35L15;
35L05
}\\
{\footnotesize{\it Key words and phrases}\/: %
fractional Laplacian, time-dependent damping, local existence, blow-up, lifespan
}

\section{Introduction}
\ \ In the present paper, we study the Cauchy problem of the $\sigma$-evolution equation with a dissipative term with a time-dependent coefficient and a power nonlinearity on the Euclidean space
\begin{equation}\label{main}
\begin{cases}
\partial_t^2u(t,x)+(-\Delta)^{\sigma/2} u(t,x)+b(t) \partial_t u(t,x) =f(u(t,x)),\ &(t,x)\ \in [0,T)\times \R^N,\\
u(0,x)=u_0(x),\ \partial_tu(0,x)=u_1(x),\ &x\ \in\ \R^N.
\end{cases}
\end{equation}
Here $N\in \N$ is the spatial dimension, $\sigma>0$ stands for the strength of diffusion and $(-\Delta)^{\sigma/2}$ is the fractional Laplacian on $\R^N$ defined by
\[
(-\Delta)^{\sigma/2} f(x):=\mathcal{F}^{-1}[|\xi|^\sigma \mathcal{F}[f](\xi)](x),
\]
where $\mathcal{F}$ and $\mathcal{F}^{-1}$ are the Fourier transform and its inverse respectively. Moreover $T\in (0,\infty]$ is the existence time of the function $u$, $u:[0,T)\times \R^N \rightarrow\R$ is an unknown function. Here $b:[0,\infty)\rightarrow \mathbb{R}_{>0}$ is a given positive function, satisfying
\begin{equation}\label{b}
    b_1(1+t)^{-\beta}\leq b(t)\leq b_2(1+t)^{-\beta},\quad \quad \mbox{and}\quad \quad |b'(t)|\leq b_3(1+t)^{-1}b(t),
\end{equation}
for $\beta\in \R$ and some constants $b_1, b_2, b_3>0$, and $f:\mathbb{R}\rightarrow\mathbb{R}$ is a nonlinear function satisfying $f(0)=0$ and there exist $p\ge 1$ and $C_0$ such that for any $z,w\in \mathbb{R}$,
\begin{equation}
\label{assumption on nonlinearity}
|f(z)-f(w)|\le C_0(|z|^{p-1}+|w|^{p-1})|z-w|.
\end{equation}
And $u_0:\mathbb{R}^N\rightarrow \mathbb{R}$ is the initial position, while $u_1:\mathbb{R}^N\rightarrow \mathbb{R}$ is the initial velocity.

Let $\mathcal{S}\left(\mathbb{R}^{N}\right)$ be the
rapidly decaying function space. For $f\in \mathcal{S}(\mathbb{R}^N)$, we define the Fourier transform of $f$ as 
\begin{equation*}
\mathcal{F}\left[ f\right] \left( \xi \right) =\widehat{f}\left( \xi \right) :=
(2\pi)^{-N/2}\int_{\R^N}f\left( x\right) e^{-ix\cdot\xi} dx,
\end{equation*}
and the inverse Fourier transform of $f$ as
\begin{equation*}
  \mathcal{F}^{-1}\left[ f\right] \left( x\right) :=(2\pi)^{-N/2}\int_{\R^N}f\left( \xi \right) e^{ix\cdot\xi } d\xi,
\end{equation*}
and extend them to $\mathcal{S}^{\prime }(\mathbb{R}^N)$ by duality. Next, we give the notations of function spaces.  For $s\in \R$ and $m\ge 0$, we denote the Sobolev space $H^{s,m}(\R^N)$ with the norm
\begin{equation*}
\|f\|_{H^{s,m}}:= \left(\int_{\R^N}\la x \ra^{2m}(\langle \nabla \rangle^s f)^2dx\right)^{1/2},
\end{equation*}
where $\langle \cdot \rangle :=(1+\vert \cdot \vert ^2)^{1/2}$.

In recent decades, the long-time behavior of the solution to the Cauchy problem related to (\ref{main}) has been investigated by many mathematicians. When $\sigma=2$ and $b(t)=1$, \eqref{main} is reduced to the semilinear damped wave equation:
\begin{equation}\label{DW}
\left\{
\begin{aligned}
\partial_t^2u(x,t)-\Delta u(x,t)+\partial_t u(x,t)=|u(x,t)|^p, &\ t>0, x\in\R^N,\\
u(x,0)=\varepsilon u_0(x),\ \partial_tu(x,0)=\varepsilon u_1(x),\ &x\ \in\ \R^N.
\end{aligned}
\right.
\end{equation}
G. Todorova and B. Yordanov \cite {TY1, TY2} studied the blow-up and
global existence of the solution to the Cauchy problem (\ref{DW}), and obtained the critical exponent $ p_c=1+\frac 2N$, which is the same as the Fujita exponent for the heat equation. Precisely, they proved the global existence for small initial data in the supercritical case
$\rho_F(N)<\rho< \frac {N+2}{[N-2]_+}$ and blow up for the data positive on average in the subcritical case $1<\rho <\rho_F(N)$.
And in the critical case $\rho=\rho_F(N)$, Qi S. Zhang \cite{Z} showed that the
solution blows up within a finite time. Moreover, the estimates of the lifespan $T(\varepsilon)$ for the blow-up solution in the subcritical and critical cases were derived by T. Li and Y. Zhou \cite{LZ}, K. Nishihara \cite{N2003}, M. Ikeda and T. Ogawa \cite{IO} and N. Lai and Y. Zhou \cite{LaiZhou} as follows:
\begin{equation*}
T(\varepsilon)\approx \left\{
\begin{array}{ll}
C \varepsilon^{-(\frac 1{p-1}-\frac N2)^{-1}},\quad &\mbox{if}\ \ 1<p<1+\frac 2N,\\
\exp(C\varepsilon^{-(p-1)}),\quad &\mbox{if}\ \ p=1+\frac 2N.
\end{array}
\right.
\end{equation*}

For the case that $\sigma=2$ and $b(t)=(1+t)^{-\beta}$, K. Nishihara \cite{N} and J. Lin, K. Nishihara and J. Zhai \cite{NLZ} investigated the Cauchy problem for the semilinear wave equation with time-dependent damping
\begin{equation}\label{1.3}
\left\{
\begin{array}{ll}
\partial_t^2u(x,t)-\Delta u(x,t)+b(t)\partial_t u(x,t)=|u(x,t)|^p, &\ t>0, x\in\R^N,\\
u(x,0)=\varepsilon u_0(x),\ \partial_tu(x,0)=\varepsilon u_1(x),\ &x\ \in\ \R^N,\\
\end{array}
\right.
\end{equation}
with $-1<\beta<1$, and they found that the critical exponent for \eqref{1.3} remains $ p_c=1+\frac 2N$. And the sharp estimate of the lifespan in the subcritical case $1<p<p_c$ was obtained by M. Ikeda and Y. Wakasugi \cite{IW}, K. Fujiwara, M. Ikeda and Y. Wakasugi \cite{FIW2019}
\begin{equation*}
T(\varepsilon)\approx C \varepsilon^{-[(\frac 1{p-1}-\frac N2)(1+\beta)]^{-1}}.
\end{equation*}
In the critical case $p=p_c$, the lower and upper estimates of the lifespan were obtained by M. Ikeda and Y. Wakasugi \cite{IW}, M. Ikeda, T. Ogawa \cite{IO} and M. Ikeda,  M. Sobajima and Y. Wakasugi\cite{ISW} as follows:
\begin{equation*}
T(\varepsilon)\approx \exp(C\varepsilon^{-(p-1)}).
\end{equation*}
Moreover, in the case $\beta = -1$, similar results were derived by K. Fujiwara, M. Ikeda and Y. Wakasugi \cite{FIW2019} and M. Ikeda and T. Inui\cite{II}:
\begin{equation*}
T(\varepsilon)\approx
\left\{
\begin{array}{ll}
C \exp{(\varepsilon^{-(\frac 1{p-1}-\frac N2)^{-1}})},\quad &\mbox{if}\ \ 1<p<1+\frac 2N,\\
\exp(\exp(C\varepsilon^{-(p-1)})),\quad &\mbox{if}\ \ p=1+\frac 2N.
\end{array}
\right.
\end{equation*}
For general effective damping, see M. D'Abbicco, S. Lucente and M. Reissig \cite{DLR}, M. D'Abbicco and S. Lucente \cite{DL}, M. Ikeda,  M. Sobajima and Y. Wakasugi\cite{ISW}, etc.

For the equation with fractional Laplacian and damping
\begin{equation}\label{1.4}
\left\{
\begin{array}{ll}
\partial_t^2u(x,t)+(-\Delta)^{\sigma /2} u(x,t)+\partial_t u(x,t)=|u(x,t)|^p, &\ t>0, x\in\R^N,\\
u(x,0)=\varepsilon u_0(x),\ \partial_tu(x,0)=\varepsilon u_1(x),\ &x\ \in\ \R^N,\\
\end{array}
\right.
\end{equation}
T. Dao and M. Reissig \cite{DR} proved that the solution would blow up in $C([0,\infty), H^\sigma(\R^N))$ for $\sigma \geq 2$ in the subcritical case $1<p<1+\frac \sigma N$ and K. Fujiwara, M. Ikeda and Y. Wakasugi \cite{FIW2020} obtained small data global well-posedness in a weighted Sobolev space $H^s(\R^N)\cap H^{0,\alpha}(\R^N)$ in the case $p>1+\frac \sigma N$ with $\alpha \geq 0$ and $\sigma>0$. Combining these two results, we could conclude that the critical exponent for the Cauchy problem \eqref{1.4} is $ p_c=1+\frac \sigma N.$ For the weakly coupled system, see \cite{DP}. Recently, similar results were improved to the case of time-dependent damping and fractional Laplacian with $\sigma >2$, by H.S. Aslan and T.A. Dao \cite{AD}. They proved the global well-posedness of small data solutions for some supercritical exponent, and obtained the upper estimate of the lifespan in subcritical case for the fractional number $\sigma>2$, and in critical case for the integer number $\sigma>2$. Moreover, with the lower bound of lifespan they got, the sharpness of the estimate of the lifespan was shown. Still the estimate of the lifespan for the case of the fractional Laplacian with $0<\sigma<2$ is open.

In this paper, we give the upper estimate of the lifespan for fractional Laplacian with $0<\sigma<2$ and time-dependent damping $\beta \in [-1, 1)$ in the subcritical and critical cases $1<p\leq p_c:=1+\frac{\sigma}{N}$. Comparing with the global existence in the supercritical case $1<p<\infty$ for $0<N\le \sigma$ and $1<p<\frac N{N-\sigma}$ for $N>\sigma$ obtained in \cite{LI}, we could conclude that $p_c=1+\frac{\sigma}{N}$ is the critical exponent for the Cauchy problem (\ref{main}) with $0<\sigma<2$. And the upper estimate of the lifespan is sharp, comparing with the lower estimate derived in \cite{LI}. Our result of the critical case is completely new even in the classical case $b(t)=1$. Besides, we obtain the upper estimate of the lifespan in the case $\beta =1$.

Before stating our local existence and lifespan estimates, let us introduce the definitions of mild and strong solutions to the Cauchy problem (\ref{main}).

\begin{definition}[Mild solution]
    Let the initial data $(u_0,u_1)\in H^{\sigma/2,m}\times H^{0,m}$, with $m\ge 0$. The function $u\in C([0,T);H^{\sigma/2,m}(\R^N))$ with $\partial_t u\in C(0,T;H^{0,m}(\R^N))$ is called the mild solution to (\ref{main}) if $u$ satisfies the following integral equation
    \begin{equation}\label{mildsol}
    u(t) = \widetilde{\mathcal{S}}_\sigma(t) u_0 + \mathcal{S}_\sigma(t) u_1 - \int_0^t \mathcal{S}_\sigma(t-\tau) (b(\tau)-1) \partial_\tau u(\tau) \, d\tau + \int_0^t \mathcal{S}_\sigma(t-\tau) |u(\tau)|^p \, d\tau.
    \end{equation}
    Here, the operator $\mathcal{S}_\sigma(t)$ is defined by
    \[
    \mathcal{S}_\sigma(t) := e^{-\frac{t}{2}}\mathcal{F}^{-1} L_\sigma(t,\xi) \mathcal{F},
    \]
    where the multiplier $L_\sigma : [0,\infty) \times \mathbb{R}^n \to \mathbb{R}$ is given by
    \[
    L_\sigma(t,\xi) :=
    \begin{cases}
    \dfrac{\sinh\left(t\sqrt{\frac{1}{4} - |\xi|^\sigma}\right)}{\sqrt{\frac{1}{4} - |\xi|^\sigma}}, & \text{if } |\xi| < 2^{-\frac{2}{\sigma}}, \\[1em]
    \dfrac{\sin\left(t\sqrt{|\xi|^\sigma - \frac{1}{4}}\right)}{\sqrt{|\xi|^\sigma - \frac{1}{4}}}, & \text{if } |\xi| > 2^{-\frac{2}{\sigma}},
    \end{cases}
    \]
    and the second linear propagator $\widetilde{\mathcal{S}}_\sigma(t)$ is defined by
    \[
    \widetilde{\mathcal{S}}_\sigma(t) := (\partial_t + 1)\mathcal{S}_\sigma(t).
    \]
\end{definition}

With the notion of mild solutions established, we now give the definition of strong solution and its lifespan, which will be used in our subsequent blow-up analysis.
\begin{definition}[Strong solution and its lifespan]
Let $N\in \mathbb{N}$, $\sigma\in (0,2)$, $(u_0,u_1)=\varepsilon (a_0, a_1)\in H^{\sigma,m}(\R^N)\times H^{\sigma/2,m}(\R^N)$ with $m\ge 0$ and $T>0$. We say that $u$ is a strong solution on $[0,T)$ if $u\in C^2([0,T);H^{0,m}(\R^N))\cap C^1([0,T);H^{\sigma/2,m}(\R^N))\cap C([0,T);H^{\sigma,m}(\R^N))$ satisfies the initial condition $u(0)=u_0$ and $\partial_tu(0)= u_1$ and satisfies the equation
\[
\partial_t^2u+(-\Delta)^{\sigma/2}u+b(t)\partial_tu=f(u)
\]
in the sense of $C([0,T);H^{0,m}(\R^N))$. Moreover, the lifespan $T(\varepsilon)=T(\varepsilon a_0, \varepsilon a_1)$ of the strong solution is defined by
\[
T(\varepsilon):=\sup\left\{\tau\in[0,T]\ \Big|\ 
u \text{ is a strong solution on }[0,\tau]
\right\}.
\]
\end{definition}

Now we are in a position to present the local existence and uniqueness theorem.
\begin{theorem}[Local existence and uniqueness]\label{thm:local}
Let $N\in \N$, $\sigma\in (0,2)$, $1<p<\infty$ for $0<N\le \sigma$ and $1<p<\frac N{N-\sigma}$ for $N>\sigma$ . Assume that $b(t)\in L^\infty(0,T)$. Then there exists a unique mild solution $u\in C([0,T),H^{\sigma/ 2,m}(\R^N))\cap C^1([0,T),H^{0,m}(\R^N))$ for some time $T$, if the initial data $(u_0, u_1)\in H^{\sigma/2,m}(\R^N)\times H^{0,m}(\R^N)$ with $m\ge 0$. Moreover, if in addition
\[
u_0\in H^{\sigma,m}(\mathbb R^N),\qquad
u_1\in H^{\sigma/2,m}(\mathbb R^N),\qquad
b\in C^1(0,T),
\]
then the above mild solution satisfies
\[
u\in C([0,T),H^{\sigma,m}(\mathbb R^N)),\qquad
\partial_t u\in C([0,T),H^{\sigma/2,m}(\mathbb R^N)),\qquad
\partial_t^2 u\in C([0,T),H^{0,m}(\mathbb R^N)),
\]
and equation \eqref{main} holds in $C([0,T),H^{0,m}(\mathbb R^N))$. 
Hence $u$ is a strong solution on $[0,T)$.
\end{theorem}

\begin{remark}
    Note that the local existence and uniqueness theorem does not require the assumption \eqref{b}.
\end{remark}

In what follows, we would consider the case of the source nonlinearity, that is,
\begin{equation*}f(u)\geq C_1 |u|^p
\end{equation*} for the exponent $p>1$ and the constant $C_1>0$.
For $N<q<N+p\sigma$, define
\begin{equation}\label{A}
A(N,p,\sigma,q):=\left(\frac 2{C_1}\right)^{p'-1}(p')^{-\frac 1p}p^{\frac {1-p'}p}\|(-\Delta)^{\frac \sigma 2}\langle x\rangle^{-q}\langle x\rangle^{\frac qp}\|^{p'/p}_{L^{p'}}\|\langle x\rangle^{-q}\|_{L^1}^{\frac 1{p'}},
\end{equation}
and
\begin{equation}\label{R}
R(\varepsilon)=A(N,p,\sigma,q)^{-\frac{p-1}{N(p-p_c)}}\left(\frac \varepsilon 4 I_0\right)^{\frac{p-1}{N(p-p_c)}}
\end{equation}
with $p_c=1+\frac \sigma N$ ,$p'=\frac p{p-1}$, and $I_0$ defined in (\ref{I0I1}) below. Besides, we assume that $b(t)$ satisfies the condition \eqref{b}.

Let us introduce the main result on the blow-up of the solution in the subcritical and critical cases.

\begin{theorem}[Blow up in the subcritical case]
Let $-1\leq \beta\leq 1$, $m\ge 0$ and $1<p<p_c$. Let $u$ be the strong solution to the Cauchy problem (\ref{main}) with the nonlinearity $f(u)$ satisfying $f(u)\geq C_1|u|^p$ for some constants $C_1>0$. Assume that the initial data $(u_0, u_1)=\varepsilon (a_0, a_1)$ with $a_0\in H^{\sigma,m}(\R^N)\cap L^1(\R^N)$ and $a_1 \in  H^{\sigma/2,m}(\R^N)\cap L^1(\R^N)$, satisfy
\begin{equation}\label{I0I1}
I_0:=\int_{\R^N}a_0(x)dx>0,\quad \quad I_1:=\int_{\R^N}a_1(x)dx>0.
\end{equation}
Let $\varepsilon_0$ be a constant small enough such that for any $0<\varepsilon\leq \varepsilon_0$
\begin{equation}\label{BLC1}
    \int_{\R^N}\langle \frac x {R(\varepsilon)}\rangle^{-q}a_0(x)dx\geq \frac 12 I_0,
\end{equation}
\begin{equation}\label{BLC2}
    \int_{\R^N}\langle \frac x {R(\varepsilon)}\rangle^{-q}a_1(x)dx\geq \frac 12 I_1,
\end{equation}
\begin{equation}\label{BLC3}
    \varepsilon I_0\leq \left(\frac {C_1} 2\right)^{-\frac 1{p-1}}\|\langle x\rangle^{-q}\|_{L^1(\R^N)}R(\varepsilon)^N.
\end{equation}
Then, for $0<\varepsilon\leq \varepsilon_0$, the lifespan of the solution can be estimated by
\begin{equation}
T(\varepsilon) \leq
\begin{cases}
C \varepsilon^{-\frac 1{\left(\frac 1{p-1}-\frac N\sigma\right)(1+\beta)}}, &\quad \mbox{for}\quad -1<\beta\leq 1,\\
\exp{\left(C\varepsilon^{-\frac{1}{\frac 1{p-1}-\frac N\sigma}}\right)}, &\quad \mbox{for}\quad \beta= -1,
\end{cases}
\end{equation}
for some constant $C$ independent of $\varepsilon$.
\end{theorem}

\begin{remark}
    Noting that for $1<p<p_c$, $R(\varepsilon)\rightarrow \infty$ as $\varepsilon\rightarrow 0$, since $\frac{p-1}{N(p-p_c)}<0$. Thus, we have
    \begin{equation*}
    \begin{aligned}
        \int_{\R^N}\langle \frac x {R(\varepsilon)}\rangle^{-q}a_0(x)dx\rightarrow \int_{\R^N}a_0(x)dx=I_0,\\
        \int_{\R^N}\langle \frac x{R(\varepsilon)}\rangle^{-q} a_1(x)dx\rightarrow \int_{\R^N}a_1(x)dx=I_1,
    \end{aligned}
    \end{equation*}
    which implies the conditions (\ref{BLC1}) and (\ref{BLC2}) are satisfied. Moreover, it follows from
    \begin{equation*}
     \varepsilon R(\varepsilon)^{-N}=A(N,p,\sigma,q)^{\frac{p-1}{(p-p_c)}}\left(\frac 14 I_0\right)^{\frac{p-1}{p_c-p}}\varepsilon^{1+\frac{p-1}{p_c-p}}
    \end{equation*}
    with the power $1+\frac{p-1}{p_c-p}>0$ that there exists a constant $\varepsilon_0$ small enough such that the condition (\ref{BLC3}) holds.
\end{remark}

\begin{theorem}[Blow up in the critical case]
Let $-1\leq \beta< 1$, $m\ge 0$ and $p= 1+\frac \sigma N$ with $0<\sigma\leq 2$. Assume that $u$ is a strong solution to the Cauchy problem \eqref{main} with the initial data $(u_0, u_1)=\varepsilon (a_0, a_1)\in \left(H^{\sigma,m}(\R^N)\cap L^1(\R^N)\right)\times \left(H^{\sigma/2,m}(\R^N)\cap L^1(\R^N)\right)$ satisfying
\begin{equation*}
\int_{\R^N}(a_0(x)+B_0 a_1(x))dx>0,\quad B_0:=\int_0^\infty\exp \left(-\int_0^t b(s)ds\right)dt.
\end{equation*}
Then, the solution $u$ blows up in a finite time. Moreover, the lifespan of the solution can be estimated by
\begin{equation}
    T(\varepsilon) \leq
    \begin{cases}
        \exp{(C\varepsilon^{-(p-1)}},\quad &\mbox{if}\quad \beta\in (-1,1),\\
        \exp{(\exp{(C\varepsilon^{-(p-1)}})},\quad &\mbox{if}\quad \beta=-1,
    \end{cases}
\end{equation}
for some constant $C$ independent of $\varepsilon$.
\end{theorem}

\begin{remark}
    Note that the coefficient $b(t)$ of the damping has no influence on the estimate of the lifespan.
\end{remark}

We briefly give the idea of the proofs. To get the local existence, we write the equation as
\begin{equation*}
\partial_t^2u+\partial_t u+(-\Delta)^{\sigma/2} u=(1-b(t)) \partial_t u +f(u),
\end{equation*}
to use the properties of the solution operators $S_\sigma(t)$ and $\tilde{S}_\sigma(t)$ for the linear damped fractional wave equation. Moreover, to overcome the difficulty of handling spatial weights coupled with the nonlinear term $\| |u|^p\|_{H^{0,m}}$, we rely on the weighted estimates for the Bessel potential operator (Lemma \ref{lem:weighted_Js}) and Sobolev inequalities. Then, based on the estimates of the operators $S_\sigma(t)$ and $\tilde{S}_\sigma(t)$ (Proposition \ref{prop:linear-bounds}), we could obtain that the mild solution exists locally by the contraction mapping principle. Furthermore, under higher-regularity conditions for the initial data and the damping coefficient, this mild solution can be further promoted to a strong solution.

To prove the blow-up of the solution, we apply the framework of the ordinary differential inequality. This method was developed by T. Li and Y. Zhou \cite{LZ}, requiring the positivity of the fundamental solution. However, this positivity condition does not hold for the damped wave equation in higher dimension case. To overcome the above difficulty, a hybrid version of the method by T. Li and Y. Zhou \cite{LZ} and the test function method by Q. Zhang \cite{Z} was introduced, which was applied to the nonlinear Schr\"{o}dinger equation by K. Fujiwara and T. Ozawa \cite{FO1,FO2}, and to the damped wave equation by K. Fujiwara, the second author and Y. Wakasugi \cite{FIW2019}. In this paper, we adapt this method to the fractional wave equation with time-dependent damping in the subcritical case. First, we construct a weighted average of the solution, 
and derive an ordinary differential inequality from the fractional nonlinear damped wave equation. Then, the weighted average of the solution is proved positive by the initial data condition and the ordinary differential inequality derived from the damped wave equation. Thus, by the blow-up and upper estimate of the lifespan of the solution to the ordinary differential inequality, we prove that the solution $u$ to the Cauchy problem \eqref{main} blows up in a finite time and obtain the upper bound of the lifespan. However, this method fails in the critical case, since the weight function make no sense if $p=1+\frac{\sigma}{N}$. Instead, we employ the idea of methods in \cite{LZ} and \cite{II} to find a new suitable weight function. Note that the damped wave equation is expected to have the diffusion phenomenon, that is, its solution behaves like that to the corresponding heat equation as $t\rightarrow \infty$, it's reasonable to use the fundamental solution to the fractional heat equation as the weight. To overcome the difficulties which come from the time-dependent damping, we multiply the equation by a nonnegative function $g(t)\in C^1([0,\infty))$ to the divergence form:
\begin{equation*}
 (g(t)u)_{tt} +(-\Delta)^{\sigma/2} (g(t)u) -(g'(t)u)_t + u_t = g(t)|u|^{\rho}.
\end{equation*}
Then, we treat the above equation as the semilinear fractional heat equation with the nonlinearity $F(x,t):=g(t)|u|^{\rho}-(g(t)u)_{tt}+(g'(t)u)_t$, and write it as
\begin{equation*}
u(t)=\Phi(G(t))\ast u_0 +\int_0^t \Phi(G(t)-G(s))\ast F(s) ds,
\end{equation*}
where $\Phi(t)=\Phi(x,t):=\mathcal{F}^{-1}[e^{-t|\xi|^\sigma}]$ is the fundamental solution to the fractional heat equation and $G(t):=\int_0^t g(s)ds$.  With the property of $\Phi(G(t))$, we get ordinary differential inequalities. Applying Lemma A.1 and Lemma A.2 in \cite{II}, we obtain the upper bound of the lifespan, which is independent of $\beta$ in the damping.

This paper is organized as follows. In Section 2, we establish some preliminaries, focusing on a fundamental weighted estimate for the Bessel potential operator. In Section 3, we apply the contraction mapping principle to prove the local existence and uniqueness of the solution. Section 4 is devoted to the subcritical case, where we use the test function method to prove the blow-up of the solution and obtain the upper bound of the lifespan. In Section 5, we investigate the critical case by rewriting the equation into a divergence form to establish the blow-up result and the upper estimate of the lifespan. 


\section{Preliminaries}

In this section, we establish a fundamental weighted estimate for the Bessel potential operator $\la \nabla \ra^{s}$ with $0<s<2$. 
The result will be used later to prove a weighted estimate of the nonlinearity.

\begin{lemma}\label{lem:weighted_Js}
Let $0<s<2$ and $a\ge 0$. Then for every function $u$ for which the right‑hand sides are finite,
\begin{align}
\norm{\la x\ra^a u}_{L^2} &\le C \norm{\la x\ra^a \la \nabla \ra^{s} u}_{L^2}, \label{eq:Jemb}\\
\norm{\la \nabla \ra^{s}(\la x\ra^a u)}_{L^2} &\le C \norm{\la x\ra^a \la \nabla \ra^{s} u}_{L^2}. \label{eq:Jcomm} 
\end{align}
The constant $C$ depends only on $N,s,a$.
\end{lemma}

\begin{proof}
We first establish \eqref{eq:Jemb}. Set $h := \la \nabla \ra^{s} u$, so that $u = \la \nabla \ra^{-s}h$. The operator $\la \nabla \ra^{-s}$ acts by convolution with the Bessel kernel $G_s$:
\[
u(x) = \int_{\R^N} G_s(x-y) h(y)\,dy.
\]
It follows from Proposition 7.6 in \cite{Hao} that  
\[
\int_{\R^N} |G_s(z)|\la z\ra^a \,dz < \infty \qquad\text{for every } a\ge 0.
\]

We use the standard estimate for moderate weights
\[
\la x\ra \le \sqrt{2}\,\la x-y\ra\,\la y\ra,
\]
which implies 
\[
\la x\ra^a \le C_a \la x-y\ra^a \la y\ra^a,\quad\text{for}\ a \ge 0.
\]

Hence, 
\begin{align*}
\la x\ra^a|u(x)|
&\le \int_{\R^N} |G_s(x-y)|\la x\ra^a |h(y)|\,dy \\
&\le C_a \int_{\R^N} |G_s(x-y)| \la x-y\ra^a \la y\ra^a |h(y)|\,dy.
\end{align*}
Applying Young's inequality, we have
\[
\|\la x\ra^a u\|_{L^2} \le C_a \|\la z\ra^a G_s(z)\|_{L^1} \|\la x\ra^a h\|_{L^2}
\le C \|\la x\ra^a \la \nabla \ra^{s} u\|_{L^2},
\]
which proves \eqref{eq:Jemb}.

We now turn to \eqref{eq:Jcomm}. Write
\[
\la \nabla \ra^{s}(\la x\ra^a u) = \la x\ra^a \la \nabla \ra^{s} u + [\la \nabla \ra^{s}, \la x\ra^a]u.
\]
It suffices to estimate the commutator. We will show
\begin{equation}\label{eq:commtarget}
 \norm{[\la \nabla \ra^{s},\la x\ra^a]\la \nabla \ra^{-s}h}_{L^2}\le C\norm{\la x\ra^a h}_{L^2},\qquad \mbox{for}\ \ h=\la \nabla \ra^s u.
\end{equation}
Set \(\alpha=s/2\in(0,1)\) and \(c_s=\alpha/\Gamma(1-\alpha)\). The semigroup representation
\begin{equation}\label{eq:semigroup}
 \la \nabla \ra^sf=f+c_s\int_0^\infty e^{-t}t^{-1-\alpha}(f-e^{t\Delta}f)\,dt.
\end{equation}
gives
\begin{equation}\label{eq:commheat}
 [\la \nabla \ra^s,\la x\ra^a]u(x)=c_s\int_0^\infty e^{-t}t^{-1-\alpha}
 \int_{\R^N}H_t(z)[\la x\ra^a-\la x-z\ra^a]u(x-z)\,dz\,dt
\end{equation}
with \(z=x-y\). Applying Taylor's theorem, we get
\begin{equation}\label{eq:taylor}
 \la x\ra^a-\la x-z\ra^a=\nabla \la x\ra^a\cdot z-Q(x,z),
\end{equation}
where
\[
 Q(x,z)=\int_0^1(1-\theta)z^TD^2\la x-\theta z\ra^a z\,d\theta.
\]
By the derivative estimates for \(\la x\ra^a\),
\begin{equation}\label{eq:Qbound}
 |Q(x,z)|\le C_a|z|^2\la z \ra^a \la x-z\ra^a.
\end{equation}
Thus \eqref{eq:commheat} can be written as \(A_s u - B_s u\), where
\[
A_s u(x) = c_s\int_0^\infty e^{-t}t^{-1-\alpha}
\int_{\R^N} H_t(z)\,\nabla \la x\ra^a\cdot z\, u(x-z)\,dz\,dt,
\]
and
\[
B_s u(x) = c_s\int_0^\infty e^{-t}t^{-1-\alpha}
\int_{\R^N} H_t(z)\, Q(x,z)\, u(x-z)\,dz\,dt.
\]

We first bound the term \(A_su\). For each \(k\), since
\[
z_kH_t(z)=-2t\,\partial_{z_k}H_t(z),
\]
integrating by parts gives
\[
\begin{aligned}
\int_{\R^N}z_kH_t(z)u(x-z)\,dz
&=-2t\int_{\R^N}\partial_{z_k}H_t(z)\,u(x-z)\,dz\\
&=2t\int_{\R^N}H_t(z)\,\partial_{z_k}u(x-z)\,dz\\
&=-2t\int_{\R^N}H_t(z)\partial_ku(x-z)\,dz\\
&=-2t\,e^{t\Delta}\partial_ku(x).
\end{aligned}
\]
Substituting this back into the expression for \(A_s u\) yields
\[
A_s u(x)
=-2c_s\sum_{k=1}^N \partial_k \la x\ra^a
\int_0^\infty e^{-t}t^{-\alpha}e^{t\Delta}\partial_ku(x)\,dt.
\]
Using the heat semigroup identity
\[
\int_0^\infty e^{-t}t^{-\alpha}e^{t\Delta}\,dt
=\Gamma(1-\alpha)\la \nabla \ra^{s-2},
\]
and the relation \(2c_s\Gamma(1-\alpha)=s\), we arrive at
\begin{equation}\label{eq:A}
A_s u
=-s\sum_{k=1}^N(\partial_k\la x\ra^a)\partial_k\la \nabla \ra^{s-2}u.
\end{equation}

Since \(|\partial_k\la x\ra^a|\le C_a\la x\ra^a\), we have
\[
\begin{aligned}
|A_s\la\nabla\ra^{-s}h(x)|
&\le C_a \la x\ra^a\sum_{k=1}^N
|\partial_k\la\nabla\ra^{-2}h(x)|\\
&\le C_a \sum_{k=1}^N\int_{\R^N}\la x\ra^a|\partial_kG_2(x-y)|\,|h(y)|\,dy\\
&\le C_a \sum_{k=1}^N \int_{\R^N}
|\partial_kG_2(x-y)|\la x-y\ra^a\la y\ra^a|h(y)|\,dy.
\end{aligned}
\]
Then, by the Bessel potential estimate and Young's inequality, we obtain
\begin{equation}\label{eq:Aest}
\|A_s\la\nabla\ra^{-s}h\|_{L^2} \le \|\la z \ra^a \partial_k G_2\|_{L^1} \|\la x\ra^a h\|_{L^2}
\le C\|\la x\ra^a h\|_{L^2}.
\end{equation}

We now turn to the remainder \(B_su\). By \eqref{eq:Qbound},
\[
 |B_s\la \nabla \ra^{-s}h(x)|
 \le C\int_0^\infty e^{-t}t^{-1-s/2}
 \int_{\R^N} H_t(z)|z|^2\la z \ra^a \la x-z\ra^a|u(x-z)|\,dz\,dt.
\]
The potential representation \(u=G_s*h\), followed by the moderate-weight inequality, gives
\[
 \la x-z\ra^a|u(x-z)|\le C\int_{\R^N} G_s(v)\la v \ra^a \la x-z-v\ra^a|h(x-z-v)|\,dv.
\]
Define the nonnegative convolution kernel
\begin{equation}\label{eq:M}
 M_{s,a}(\zeta):=C\int_0^\infty e^{-t}t^{-1-s/2}
 \int_{\R^N}H_t(z)|z|^2\la z \ra^a
 G_s(\zeta-z)\la \zeta-z \ra^{a}\,dz\,dt.
\end{equation}
Then
\begin{equation}\label{eq:Bpoint}
 |B_s\la \nabla \ra^{-s}h(x)|\le (M_{s,a}*(\la\cdot\ra^a|h|))(x).
\end{equation}
Tonelli's theorem and the Gaussian kernel estimate
\[
\int_{\R^N} H_t(z)|z|^2\la z \ra^a\,dz\le C(t+t^{1+a/2})
\]
show
\[
 \norm{M_{s,a}}_{L^1}
 \le C\int_0^\infty e^{-t}t^{-1-s/2}
 \big(t+t^{1+a/2}\big)\,dt
 =C\left[\Gamma(1-s/2)+\Gamma(1+(a-s)/2)\right]<\infty.
\]
Both Gamma arguments are positive because \(0<s<2\) and \(a\ge0\). Hence Young's inequality
applied to \eqref{eq:Bpoint} yields
\begin{equation}\label{eq:Best}
 \norm{B_s\la \nabla \ra^{-s}h}_{L^2}\le C\norm{\la x\ra^a h}_{L^2}.
\end{equation}
Combining \eqref{eq:Aest} and \eqref{eq:Best} proves \eqref{eq:commtarget}.

Finally,
\[
 \norm{\la \nabla \ra^{s}(\la x\ra^a u)}_{L^2}
 \le\norm{\la x\ra^a \la \nabla \ra^{s} u}_{L^2}
 +\norm{[\la \nabla \ra^{s},\la x\ra^a]u}_{L^2}
 \le C\norm{\la x\ra^a h}_{L^2}
 =C\norm{\la x\ra^a \la \nabla \ra^{s} u}_{L^2}.
\]
This proves \eqref{eq:Jcomm}.

\end{proof}

\section{Proof of local existence and uniqueness}

In this section, we apply the contraction mapping principle to get the local existence of the solution to the Cauchy problem (\ref{main}). The proof is based on the estimates of the operators $S_\sigma(t)$ and $\tilde{S}_\sigma(t)$, which is stated as follows.

\begin{proposition}\label{prop:linear-bounds}
Let $m \ge 0$. For any $t \ge 0$, the following linear operator bounds hold:
\begin{align}
\|\mathcal S_\sigma(t) u_1\|_{H^{\frac{\sigma}{2}, m}} &\le C \langle t \rangle^{\frac{m}{\sigma}} \|u_1\|_{H^{0,m}}, \label{prop-1} \\
\|\widetilde{\mathcal{S}}_\sigma(t) u_0\|_{H^{\frac{\sigma}{2}, m}} &\le C \langle t \rangle^{\frac{m}{\sigma}} \|u_0\|_{H^{\frac{\sigma}{2}, m}}, \label{prop-2} \\
\|\partial_t \mathcal S_\sigma(t) g\|_{H^{0, m}} &\le C \langle t \rangle^{\frac{m}{\sigma}} \|g\|_{H^{0,m}}, \label{prop-3} \\
\|\mathcal S_\sigma(t) h\|_{H^{0, m}} &\le C \langle t \rangle^{\frac{m}{\sigma}} \|h\|_{H^{-\frac{\sigma}{2}, m}}. \label{prop-4}
\end{align}
\end{proposition}

\begin{proof}
We prove \eqref{prop-1} in detail, and \eqref{prop-2}--\eqref{prop-4} follow analogously. Splitting the norm into unweighted and weighted parts, we get
\[
\| \langle \cdot \rangle^m (1-\Delta)^{\frac{\sigma}{4}} \mathcal S_\sigma(t) u_1 \|_{L^2} \le \|(1-\Delta)^{\frac{\sigma}{4}} \mathcal S_\sigma(t) u_1 \|_{L^2} + \| |\cdot|^m (1-\Delta)^{\frac{\sigma}{4}} \mathcal S_\sigma(t) u_1 \|_{L^2}.
\]

For the unweighted low-frequency part $|\xi| < 2^{-\frac{2}{\sigma}}$, using the algebraic inequality 
\[
\sqrt{1/4 - |\xi|^\sigma} - 1/2 \le -|\xi|^\sigma,
\] 
we have
\[
\norm{e^{-\frac{t}{2}} L_\sigma(t,\xi)\langle\xi\rangle^{\frac{\sigma}{2}}\hat u_1}_{L^2(|\xi|<2^{-2/\sigma})} \le \norm{e^{-t|\xi|^\sigma}\langle\xi\rangle^{\frac{\sigma}{2}}}_{L^\infty(|\xi|<2^{-2/\sigma})} \norm{\hat u_1}_{L^2} \le C \norm{u_1}_{L^2}.
\]
For the unweighted high-frequency part $|\xi| > 2^{-\frac{2}{\sigma}}$, since
\[
\norm{e^{-\frac{t}{2}}\langle\xi\rangle^{\frac{\sigma}{2}} L_\sigma(t,\xi)}_{L^\infty(|\xi|>2^{-\frac{2}{\sigma}})}
\le C e^{-\frac{t}{4}},
\]
we obtain
\begin{align*}
\norm{e^{-\frac{t}{2}} L_\sigma(t,\xi)\langle\xi\rangle^{\frac{\sigma}{2}}\hat u_1}_{L^2(|\xi|>2^{-\frac{2}{\sigma}})} 
&\le \norm{e^{-\frac{t}{2}}\langle\xi\rangle^{\frac{\sigma}{2}} L_\sigma(t,\xi)}_{L^\infty(|\xi|>2^{-\frac{2}{\sigma}})} \norm{\hat u_1}_{L^2} \\
&\le C e^{-\frac{t}{4}}\norm{u_1}_{L^2}.
\end{align*}
Combining the two parts yields
\[
\norm{(1-\Delta)^{\frac{\sigma}{4}} \mathcal S_\sigma(t) u_1}_{L^2}
\le Ct \norm{u_1}_{L^2}.
\]

For the weighted part, we first estimate the low-frequency component by introducing a smooth cutoff function $\chi \in C_c^\infty(\R^N)$ such that $\chi(\xi)=1$ for $|\xi|\le 1$ and $\chi(\xi)=0$ for $|\xi|>2$. Define the low-frequency kernel
\[
K(t,x) := e^{-\frac{t}{2}} \mathcal F^{-1}\big( L_\sigma(t,\xi)\langle\xi\rangle^{\frac{\sigma}{2}}\chi(\xi) \big)(x).
\]
We claim that for any $\beta \ge 0$, $\||\cdot |^{\beta} K(t)\|_{L^{1}} \le C \langle t \rangle^{\frac{\beta}{\sigma}}$. For $0 < t < 1$, this bound follows directly from the smoothness and compact support of $L_\sigma(t,\xi)\la \xi\ra^{\frac{\sigma}{2}}\chi(\xi)$. For $t \ge 1$, we introduce the self-similar rescaling $\xi = t^{-1/\sigma} \eta$ and $x = t^{1/\sigma} y$:
\begin{align*}
\||\cdot |^{\beta}K(t)\|_{L^{1}}
&= \int_{\mathbb{R}^{N}}|x|^{\beta} \left|\int_{\mathbb{R}^{N}}e^{i x\cdot\xi}e^{-\frac{t}{2}}L_\sigma(t,\xi)\langle\xi\rangle^{\frac{\sigma}{2}}\chi (\xi)\,d\xi \right|\,dx \\
&= t^{\frac{\beta}{\sigma}} \int_{\mathbb{R}^{N}}|y|^{\beta} \left|\int_{\mathbb{R}^{N}}e^{i y\cdot\eta}e^{-\frac{t}{2}}L_\sigma(t,t^{-1/\sigma}\eta)\langle t^{-1/\sigma}\eta\rangle^{\frac{\sigma}{2}}\chi (t^{-1/\sigma}\eta)\,d\eta \right|\,dy.
\end{align*}
Using $|y|^{\beta} \le \langle y \rangle^{2n - (N + 1)}$ for a sufficiently large integer $n \in \mathbb{N}$, integration by parts yields
\[
\||\cdot |^{\beta}K(t)\|_{L^{1}} \le C t^{\frac{\beta}{\sigma}} \int_{\mathbb{R}^{N}}\langle y\rangle^{-(N + 1)} \left|\int_{\mathbb{R}^{N}}e^{iy\cdot\eta}\langle \nabla_{\eta}\rangle^{2n} \left(e^{-\frac{t}{2}}L_\sigma(t,t^{-1/\sigma}\eta)\langle t^{-1/\sigma}\eta\rangle^{\frac{\sigma}{2}}\chi (t^{-1/\sigma}\eta)\right)d\eta \right|\,dy.
\]
By the definition of $L_\sigma(t,\xi)$, for $t \ge 1$, the function satisfies the exponential bound
\[
\left|\langle \nabla_{\eta}\rangle^{2n}\left(e^{-\frac{t}{2}}L_\sigma(t,t^{-1/\sigma}\eta)\langle t^{-1/\sigma}\eta\rangle^{\frac{\sigma}{2}}\chi (t^{-1/\sigma}\eta)\right)\right|\le C e^{-c|\eta|^{\sigma}}.
\]
Integrating over $\eta$ and $y$ gives $\||\cdot |^{\beta}K(t)\|_{L^{1}} \le C t^{\frac{\beta}{\sigma}}$. Applying this bound with $\beta = m$ and $\beta = 0$, we can estimate the low-frequency weighted part using Young's convolution inequality:
\begin{equation}\label{wl}
\begin{aligned}
\norm{|\cdot|^m \mathcal F^{-1}\big( e^{-\frac{t}{2}} L_\sigma(t,\xi)\chi(\xi)\langle\xi\rangle^{\frac{\sigma}{2}}\hat u_1 \big)}_{L^2}
&= \norm{|x|^m \int_{\R^N} K(t,x-y) u_1(y)\,dy}_{L^2} \\
&\le C \norm{|\cdot|^m K(t)}_{L^1} \norm{u_1}_{L^2} + C \norm{K(t)}_{L^1} \norm{|\cdot|^m u_1}_{L^2} \\
&\le C \langle t \rangle^{\frac{m}{\sigma}} \|u_1\|_{H^{0,m}}.
\end{aligned}
\end{equation}

For high frequencies, applying Plancherel's Theorem, we get
\begin{align*}
&\norm{|\cdot|^m \mathcal F^{-1}\big( e^{-\frac{t}{2}} L_\sigma(t,\xi)(1-\chi(\xi))\langle\xi\rangle^{\frac{\sigma}{2}}\hat u_1 \big)}_{L^2}\\
&= \norm{|\nabla_\xi|^m \big[ e^{-\frac{t}{2}} L_\sigma(t,\xi)(1-\chi(\xi))\langle\xi\rangle^{\frac{\sigma}{2}}\hat u_1 \big]}_{L^2} \\
&\le C \bigg( \sum_{j=1}^n \sum_{k=0}^{[m]+1} \norm{\partial_{\xi_j}^k \big( e^{-\frac{t}{2}} L_\sigma(t,\xi)(1-\chi(\xi))\langle\xi\rangle^{\frac{\sigma}{2}} \big)}_{L^\infty} \bigg) \norm{\langle\nabla_\xi\rangle^m \hat u_1}_{L^2} \\
&\le C \norm{\langle\cdot\rangle^m u_1}_{L^2}.
\end{align*}
Combining these bounds proves \eqref{prop-1}.
\end{proof}

Now, we are ready to establish the local existence and uniqueness of the mild solution to the system \eqref{main}. For a small local time $T > 0$, define 
\[
X_T = \left\{ u \in C([0,T]; H^{\frac{\sigma}{2}, m}) \ \Big|\ \partial_t u \in C([0,T]; H^{0, m}) \right\},
\]
equipped with the norm
\[
\|u\|_{X_T} = \sup_{t \in [0,T]} \left( \|u(t)\|_{H^{\frac{\sigma}{2}, m}} + \|\partial_t u(t)\|_{H^{0, m}} \right).
\]
We define the solution operator $\mathcal{T}(u)(t)$ as the right-hand side of the mild solution formulation \eqref{mildsol}. To apply the fixed-point theorem, we first establish a priori estimates for $u(t)$ and its time derivative $\partial_t u(t)$ to show that $\mathcal{T}$ is a contraction on a suitably chosen closed ball in $X_T$.

We bound each of the four terms comprising $\mathcal{T}(u)(t)$ in $H^{\sigma/2, m}$. For the initial terms, Proposition \ref{prop:linear-bounds} gives
\[
\|\widetilde{\mathcal S}_\sigma(t) u_0\|_{H^{\frac{\sigma}{2},m}} + \|\mathcal S_\sigma(t) u_1\|_{H^{\frac{\sigma}{2},m}} \le C \langle t \rangle^{\frac{m}{\sigma}} \left( \|u_0\|_{H^{\frac{\sigma}{2},m}} + \|u_1\|_{H^{0,m}} \right).
\]
For the linear integral term, applying \eqref{prop-1}, we obtain
\begin{equation}\label{thirdle}
\begin{aligned}
\left\| \int_0^t \mathcal S_\sigma(t-\tau) (b(\tau)-1) \partial_\tau u \, d\tau \right\|_{H^{\frac{\sigma}{2}, m}}
&\le \int_0^t \| \mathcal S_\sigma(t-\tau) (b(\tau)-1) \partial_\tau u \|_{H^{\frac{\sigma}{2}, m}} d\tau \\
&\le \int_0^t |b(\tau)-1| \langle t-\tau \rangle^{\frac{m}{\sigma}} \| \partial_\tau u \|_{H^{0,m}} d\tau \\
&\le C t \left(\sup_{\tau \in [0,t]} |b(\tau)-1|\right) \langle t \rangle^{\frac{m}{\sigma}} \sup_{\tau \in [0,t]} \| \partial_t u(\tau) \|_{H^{0,m}}.
\end{aligned}
\end{equation}
For the nonlinear integral term, it follows from \eqref{prop-1} that
\begin{equation}\label{forthle}
\left\| \int_0^t \mathcal S_\sigma(t-\tau) |u|^p \, d\tau \right\|_{H^{\frac{\sigma}{2}, m}} \le \int_0^t \langle t-\tau \rangle^{\frac{m}{\sigma}} \||u|^p\|_{H^{0,m}} d\tau \le C t \langle t \rangle^{\frac{m}{\sigma}} \sup_{\tau \in [0,t]} \| \langle x \rangle^m |u|^p \|_{L^2}.
\end{equation}
To handle the spatial weight, Sobolev inequality and lemma \ref{lem:weighted_Js} yield
\begin{equation}\label{nonlinearle}
\begin{aligned}
\| \langle x \rangle^m |u|^p \|_{L^2} &= \| \langle x \rangle^{\frac{m}{p}} u \|_{L^{2p}}^p \le C \| \langle x \rangle^{\frac{m}{p}} u \|_{H^{\frac{\sigma}{2}, 0}}^{p\theta} \| \langle x \rangle^{\frac{m}{p}} u \|_{L^2}^{p(1-\theta)} \\
&\le C \|u\|_{H^{\frac{\sigma}{2}, \frac{m}{p}}}^p \le C \|u\|_{H^{\frac{\sigma}{2}, m}}^p.
\end{aligned}
\end{equation}
Inserting \eqref{nonlinearle} into \eqref{forthle} gives
\begin{equation}\label{eq:u_nonlinear_final}
\left\| \int_0^t \mathcal S_\sigma(t-\tau) |u|^p \, d\tau \right\|_{H^{\frac{\sigma}{2}, m}} \le C t \langle t \rangle^{\frac{m}{\sigma}} \sup_{\tau \in [0,t]} \|u(\tau)\|_{H^{\frac{\sigma}{2}, m}}^p.
\end{equation}

In order to control the solution in $X_T$, we also estimate the time derivative $v(t) := \partial_t u(t)$. Notice that at $t=0$, the operators $\widetilde{\mathcal S}_\sigma(t)$ and $\mathcal S_\sigma(t)$ satisfy initial conditions
\[
\mathcal S_\sigma(0)=0,\qquad \partial_t\mathcal S_\sigma(0)=I,
\]
and
\[
\widetilde{\mathcal S}_\sigma(0)=I,\qquad \partial_t\widetilde{\mathcal S}_\sigma(0)=0.
\]
Besides, the following differential identities hold
\begin{equation}\label{propagator-identities}
\partial_t\mathcal S_\sigma(t) = \widetilde{\mathcal S}_\sigma(t) - \mathcal S_\sigma(t), \qquad \partial_t\widetilde{\mathcal S}_\sigma(t) = -(-\Delta)^{\sigma/2}\mathcal S_\sigma(t).
\end{equation}
Consequently, taking second time derivatives yields
\begin{equation}\label{second-propagator-identities}
\partial_t^2\mathcal S_\sigma(t) = -(-\Delta)^{\sigma/2}\mathcal S_\sigma(t) - \partial_t\mathcal S_\sigma(t), \qquad \partial_t^2\widetilde{\mathcal S}_\sigma(t) = -(-\Delta)^{\sigma/2}\widetilde{\mathcal S}_\sigma(t) - \partial_t \widetilde{\mathcal S}_\sigma(t).
\end{equation}

Directly differentiating the mild solution formula \eqref{mildsol} yields the representation
\begin{equation}\label{veq-direct}
v(t) = -\mathcal S_\sigma(t) (-\Delta)^{\sigma/2} u_0 + \partial_t \mathcal S_\sigma(t) u_1 - \int_0^t \partial_t \mathcal S_\sigma(t-\tau) (b(\tau)-1) v(\tau) \, d\tau + \int_0^t \partial_t \mathcal S_\sigma(t-\tau) |u(\tau)|^p \, d\tau.
\end{equation}
Following the similar argument as for the estimate of $u(t)$ in \eqref{mildsol}, we apply Proposition 3.1 to obtain

\begin{equation}\label{v-final-estimate}
\begin{aligned}
\|v(t)\|_{H^{0,m}} \le &C \langle t \rangle^{\frac{m}{\sigma}} \left( \|u_0\|_{H^{\frac{\sigma}{2},m}} + \|u_1\|_{H^{0,m}} \right)\\
& + C t \langle t \rangle^{\frac{m}{\sigma}} \left( \sup_{\tau \in [0,t]} |b(\tau)-1| \sup_{\tau \in [0,t]} \|v(\tau)\|_{H^{0,m}} + \sup_{\tau \in [0,t]} \|u(\tau)\|_{H^{\frac{\sigma}{2},m}}^p \right).
\end{aligned}
\end{equation}

Combining the bounds for $u(t)$ and $\partial_t u(t)$, we obtain the following estimate for the solution operator $\mathcal{T}(u)$:
\begin{equation}\label{Phi-bound}
\begin{aligned}
\|\mathcal{T}(u)\|_{X_T} &\le C_1 \langle T \rangle^{\frac{m}{\sigma}} \left( \|u_0\|_{H^{\frac{\sigma}{2}, m}} + \|u_1\|_{H^{0, m}} \right) \\
&\quad + C_2 T \langle T \rangle^{\frac{m}{\sigma}} \left(\sup_{\tau \in [0,T]} |b(\tau)-1|\right) \|u\|_{X_T} + C_3 T \langle T \rangle^{\frac{m}{\sigma}} \|u\|_{X_T}^p.
\end{aligned}
\end{equation}

Consider the closed ball $B_R = \{ u \in X_T \mid \|u\|_{X_T} \le R \}$ with radius chosen as
\[
R = 2^{\frac{m}{2\sigma}+1} C_1 \left( \|u_0\|_{H^{\frac{\sigma}{2}, m}} + \|u_1\|_{H^{0, m}} \right).
\]
By selecting $T$ small enough such that
\[
C_2 T \langle T \rangle^{\frac{m}{\sigma}} \sup_{\tau \in [0,T]} |b(\tau)-1| \le \frac{1}{4} \quad \text{and} \quad C_3 T \langle T \rangle^{\frac{m}{\sigma}} R^{p-1} \le \frac{1}{4},
\]
we guarantee that $\|\mathcal{T}(u)\|_{X_T} \le \frac{R}{2} + \frac{R}{4} + \frac{R}{4} = R$. Hence, $\mathcal{T}(B_R) \subset B_R$.

Furthermore, for any $u, w \in B_R$, a parallel calculation for the difference $\mathcal{T}(u) - \mathcal{T}(w)$ gives:
\[
\|\mathcal{T}(u) - \mathcal{T}(w)\|_{X_T} \le \left( C_2 T \langle T \rangle^{\frac{m}{\sigma}} \sup_{\tau \in [0,T]} |b(\tau)-1| + C_4 T \langle T \rangle^{\frac{m}{\sigma}} R^{p-1} \right) \|u - w\|_{X_T}.
\]
By choosing $T > 0$ sufficiently small, we can ensure that the Lipschitz constant is strictly less than $1$, proving that $\mathcal{T}$ is a strict contraction on $B_R$. Therefore, by the fixed-point theorem, there exists a unique local mild solution $u \in X_T$ to system \eqref{main}.

Having established the local existence and uniqueness of the mild solution in \(X_T\), we now prove that it is actually a strong solution when the data and the damping 
coefficient are more regular. Assume from now on that
\[
u_0\in H^{\sigma,m},\qquad u_1\in H^{\sigma/2,m},\qquad b\in W^{1,\infty}(0,T).
\]
Let $u\in X_T$ be the mild solution constructed in the previous section, and $v(t)=\partial_t u(t)$. We shall show that $v\in C([0,T];H^{\sigma/2,m})$ and 
$\partial_t v\in C([0,T];H^{0,m})$. Once this is proved, the equation
\[
(-\Delta)^{\sigma/2}u=|u|^p-b(t)\partial_t u-\partial_t^2 u
\]
will give $u\in C([0,T];H^{\sigma,m})$.

Differentiating \eqref{mildsol} with respect to $t$ and integrating by parts, we obtain
\begin{equation}\label{veq-again}
\begin{aligned}
v(t)
&=\widetilde{\mathcal S}_\sigma(t)u_1+\mathcal S_\sigma(t)u_2\\
&\quad-\int_0^t\mathcal S_\sigma(t-\tau)(b(\tau)-1)\partial_\tau v(\tau)\,d\tau
-\int_0^t\mathcal S_\sigma(t-\tau)b'(\tau)v(\tau)\,d\tau\\
&\quad+\int_0^t\mathcal S_\sigma(t-\tau)p|u(\tau)|^{p-2}u(\tau)v(\tau)\,d\tau,
\end{aligned}
\end{equation}
where
\[
u_2:=-(-\Delta)^{\sigma/2}u_0-b(0)u_1+|u_0|^p\in H^{0,m}.
\]
The right-hand side of \eqref{veq-again} still contains $u$. 
In order to close the fixed-point argument, we reconstruct $u$ from $v$. 
Define
\begin{equation}\label{Uv-def}
U[v](t):=u_0+\int_0^t v(s)\,ds.
\end{equation}
Then $U[v]$ is determined by $v$, and if $v=\partial_t u$, then $U[v]=u$. 
Replacing $u$ by $U[v]$ in \eqref{veq-again}, we obtain the closed equation
\begin{equation}\label{v-closed}
\begin{aligned}
v(t)
&=\widetilde{\mathcal S}_\sigma(t)u_1+\mathcal S_\sigma(t)u_2\\
&\quad-\int_0^t\mathcal S_\sigma(t-\tau)(b(\tau)-1)\partial_\tau v(\tau)\,d\tau
-\int_0^t\mathcal S_\sigma(t-\tau)b'(\tau)v(\tau)\,d\tau\\
&\quad+\int_0^t\mathcal S_\sigma(t-\tau)p|U[v](\tau)|^{p-2}U[v](\tau)v(\tau)\,d\tau.
\end{aligned}
\end{equation}
Now every term on the right-hand side depends only on $v$.

For a small time $T_0$, define the space
\[
Y_{T_0}=\left\{v\in C([0,T_0];H^{\sigma/2,m})\ \Big|\ 
\partial_t v\in C([0,T_0];H^{0,m})\right\},
\]
equipped with the norm
\[
\|v\|_{Y_{T_0}}
=\sup_{t\in[0,T_0]}\left(
\|v(t)\|_{H^{\sigma/2,m}}+\|\partial_t v(t)\|_{H^{0,m}}
\right).
\]
Let $\Psi(v)$ denote the right-hand side of \eqref{v-closed}. 
We show that $\Psi$ is a strict contraction on a closed ball of $Y_{T_0}$ 
for $T_0>0$ sufficiently small.

It follows from \eqref{Uv-def} that
\[
\sup_{t\in[0,T_0]}\|U[v](t)\|_{H^{\sigma/2,m}}
\le \|u_0\|_{H^{\sigma/2,m}}+T_0\|v\|_{Y_{T_0}}.
\]
Then, applying Proposition \ref{prop:linear-bounds}, we drive the estimate for $\Psi(v)$ as follows:
\[
\|\widetilde{\mathcal S}_\sigma(t)u_1\|_{H^{\sigma/2,m}}
\le C\langle t\rangle^{m/\sigma}\|u_1\|_{H^{\sigma/2,m}},
\]
and
\[
\|\mathcal S_\sigma(t)u_2\|_{H^{\sigma/2,m}}
\le C\langle t\rangle^{m/\sigma}\|u_2\|_{H^{0,m}}.
\]
For the damping terms, using \eqref{prop-1} and $b\in C^1(0,T)$,
\[
\left\|\int_0^t\mathcal S_\sigma(t-\tau)(b(\tau)-1)\partial_\tau v(\tau)\,d\tau
\right\|_{H^{\sigma/2,m}}\le C T_0\langle T_0\rangle^{m/\sigma}
\|b-1\|_{L^\infty}\|\partial_t v\|_{C([0,T_0];H^{0,m})},
\]
and
\[
\left\|\int_0^t\mathcal S_\sigma(t-\tau)b'(\tau)v(\tau)\,d\tau
\right\|_{H^{\sigma/2,m}}\le C T_0\langle T_0\rangle^{m/\sigma}
\|b'\|_{L^\infty}\|v\|_{C([0,T_0];H^{0,m})}.
\]
For the nonlinear term, from \eqref{prop-1} and \eqref{nonlinearle}, we obtain
\[
\begin{aligned}
\left\|\int_0^t\mathcal S_\sigma(t-\tau)p|U[v](\tau)|^{p-2}U[v](\tau)v(\tau)\,d\tau \right\|_{H^{\sigma/2,m}}&\le 
C T_0\langle T_0\rangle^{m/\sigma} \|p|U[v]|^{p-2}U[v]v\|_{H^{0,m}}\\
&\le C T_0\langle T_0\rangle^{m/\sigma}\|U[v]\|_{H^{\sigma/2,m}}^{p-1}\|v\|_{H^{\sigma/2,m}}.
\end{aligned}
\]
Moreover, by \eqref{Uv-def},
\[
\sup_{t\in[0,T_0]}\|U[v](t)\|_{H^{\sigma/2,m}}
\le \|u_0\|_{H^{\sigma/2,m}}+T_0\|v\|_{Y_{T_0}}.
\]
Therefore,
\begin{equation}\label{Psi-v-estimate}
\begin{aligned}
\sup_{t\in[0,T_0]} \|\Psi(v)(t)\|_{H^{\frac{\sigma}{2},m}} \le &C \langle T_0 \rangle^{\frac{m}{\sigma}} \left( \|u_1\|_{H^{\frac{\sigma}{2},m}} + \|u_2\|_{H^{0,m}} \right)\\
& + C T_0 \langle T_0 \rangle^{\frac{m}{\sigma}} \left( \|b-1\|_{L^\infty} \sup_{t\in[0,T_0]} \|\partial_t v(t)\|_{H^{0,m}} + \|b'\|_{L^\infty} \sup_{t\in[0,T_0]} \|v(t)\|_{H^{0,m}} \right)\\
& + C T_0 \langle T_0 \rangle^{\frac{m}{\sigma}} \left( \|u_0\|_{H^{\frac{\sigma}{2},m}} + T_0 \|v\|_{Y_{T_0}} \right)^{p-1} \sup_{t\in[0,T_0]} \|v(t)\|_{H^{\frac{\sigma}{2},m}}.
\end{aligned}
\end{equation}

Similarly, differentiating the equation \eqref{v-closed} with respect to $t$ and using the initial condition $\mathcal S_\sigma(0)=0$, we have
\begin{equation}\label{dt-v-eq}
\begin{aligned}
\partial_t v(t) &= \partial_t \widetilde{\mathcal S}_\sigma(t)u_1 + \partial_t \mathcal S_\sigma(t)u_2 \\
&\quad - \int_0^t \partial_t \mathcal S_\sigma(t-\tau)(b(\tau)-1)\partial_\tau v(\tau)\,d\tau \\
&\quad - \int_0^t \partial_t \mathcal S_\sigma(t-\tau)b'(\tau)v(\tau)\,d\tau \\
&\quad + \int_0^t \partial_t \mathcal S_\sigma(t-\tau)p|U[v](\tau)|^{p-2}U[v](\tau)v(\tau)\,d\tau.
\end{aligned}
\end{equation}
Applying Proposition \ref{prop:linear-bounds}, we obtain the corresponding estimate for the time derivative in $H^{0,m}$:
\begin{equation}\label{dt-Psi-v-estimate}
\begin{aligned}
\sup_{t\in[0,T_0]} \|\partial_t \Psi(v)(t)\|_{H^{0,m}} \le &C \langle T_0 \rangle^{\frac{m}{\sigma}} \left( \|u_1\|_{H^{\frac{\sigma}{2},m}} + \|u_2\|_{H^{0,m}} \right)\\
& + C T_0 \langle T_0 \rangle^{\frac{m}{\sigma}} \left( \|b-1\|_{L^\infty} \sup_{t\in[0,T_0]} \|\partial_t v(t)\|_{H^{0,m}} + \|b'\|_{L^\infty} \sup_{t\in[0,T_0]} \|v(t)\|_{H^{0,m}} \right)\\
& + C T_0 \langle T_0 \rangle^{\frac{m}{\sigma}} \left( \|u_0\|_{H^{\frac{\sigma}{2},m}} + T_0 \|v\|_{Y_{T_0}} \right)^{p-1} \sup_{t\in[0,T_0]} \|v(t)\|_{H^{\frac{\sigma}{2},m}}.
\end{aligned}
\end{equation}

Therefore, on a ball $\|v\|_{Y_{T_0}}\le R$,
\[
\begin{aligned}
\|\Psi(v)\|_{Y_{T_0}}
&\le C_1\langle T_0\rangle^{m/\sigma}
\left(\|u_1\|_{H^{\sigma/2,m}}+\|u_2\|_{H^{0,m}}\right)\\
&\quad+C_2T_0\langle T_0\rangle^{m/\sigma}
\left(\|b-1\|_{L^\infty}+\|b'\|_{L^\infty}\right)\|v\|_{Y_{T_0}}\\
&\quad+C_3T_0\langle T_0\rangle^{m/\sigma}
\left(\|u_0\|_{H^{\sigma/2,m}}+T_0\|v\|_{Y_{T_0}}\right)^{p-1}
\|v\|_{Y_{T_0}}.
\end{aligned}
\]
Similarly, for any $v,w$ in the same ball,
\[
\|\Psi(v)-\Psi(w)\|_{Y_{T_0}}
\le C_R T_0\langle T_0\rangle^{m/\sigma}\|v-w\|_{Y_{T_0}},
\]
where $C_R$ depends on $R$, $\|b\|_{C^1}$, and $\|u_0\|_{H^{\sigma/2,m}}$, 
but not on $T_0$ when $0<T_0\le1$. 

Choose
\[
R=2^{\frac{m}{2\sigma}+1} C_1
\left(\|u_1\|_{H^{\sigma/2,m}}+\|u_2\|_{H^{0,m}}+1\right),
\]
and then take $T>0$ sufficiently small such that
\[
C_2 T_0\langle T_0\rangle^{m/\sigma}\left(\|b-1\|_{L^\infty}+\|b'\|_{L^\infty}\right)\le\frac14 \quad \text{and} \quad C_2 T_0\langle T_0\rangle^{m/\sigma}\le\frac14.
\]
Then $\Psi$ maps the closed ball $\{v\in Y_{T_0}:\|v\|_{Y_{T_0}}\le R\}$ into itself 
and is a strict contraction. By the fixed-point theorem, there exists a unique 
$v\in Y_{T_0}$ solving \eqref{v-closed}.

Now define
\[
\widetilde u(t):=U[v](t)=u_0+\int_0^t v(s)\,ds.
\]
Then $\partial_t\widetilde u=v$. Integrating the differentiated equation in time and using the definition of $u_2$ show that $\widetilde u$ satisfies the original mild formula \eqref{mildsol} with the same initial data. By uniqueness of the mild solution in $X_T$, we conclude that $\widetilde u=u$ on $[0,T_0]$. Consequently,
\[
v=\partial_t u\in C([0,T_0];H^{\sigma/2,m}),\qquad
\partial_t v=\partial_t^2 u\in C([0,T_0];H^{0,m}).
\]
Finally, since the PDE holds in $H^{0,m}$,
\[
(-\Delta)^{\sigma/2}u=|u|^p-b(t)\partial_t u-\partial_t^2 u
\in C([0,T_0];H^{0,m}),
\]
which implies
$(-\Delta)^{\sigma/2}$ yields
\[
u\in C([0,T_0];H^{\sigma,m}).
\]
Thus the mild solution is a strong solution on $[0,T_0]$. 

We now prove that the higher regularity persists on every compact subinterval of \([0,T_{\max})\). Fix an arbitrary \(T\) such that \(0<T<T_{\max}\). Let $T(\varepsilon)$ is the lifespan of the strong solution. It's known that \(T(\varepsilon)>0\). Assume, for the sake of contradiction, that \(T(\varepsilon)<T\).

Since the higher-regularity solution is continuous up to \(T(\varepsilon)\) from the left, the left limits \(u(T(\varepsilon))\) and \(\partial_t u(T(\varepsilon))\) belong to \(H^{\sigma,m}\) and \(H^{\sigma/2,m}\), respectively. Also \(b\in C^1([0,T])\). Hence we may restart the above argument at the initial time \(T_*\) with initial data 
\[
(u(T(\varepsilon)),\partial_t u(T(\varepsilon)))\in H^{\sigma,m}\times H^{\sigma/2,m}.
\]
This restart yields a higher-regularity solution on \([T(\varepsilon),T(\varepsilon)+\delta]\) for some \(\delta>0\).

By uniqueness of the mild solution in the class \(X_T\), this restarted solution coincides with the original maximal mild solution on their common interval of existence. Therefore the higher regularity extends strictly beyond \(T_*\), which contradicts the definition of \(T(\varepsilon)\) as a supremum. Hence \(T(\varepsilon)=T\).

Since \(T<T_{\max}\) was arbitrary, we conclude that
\[
u\in C([0,T_{\max});H^{\sigma,m}),\qquad
\partial_t u\in C([0,T_{\max});H^{\sigma/2,m}),\qquad
\partial_t^2 u\in C([0,T_{\max});H^{0,m}).
\]
This completes the proof of Theorem \ref{thm:local}.

\section{Proof of blow-up and upper bound of lifespan in subcritical case}

In this section, we apply the test function method to prove the blow-up of the solution in the subcritical case. We focus on deriving an ordinary differential inequality for the weighted average of the solution. Then, by the comparison lemma in \cite{FIW2019}, we could obtain the blow-up result and the upper bound of the lifespan. 

Define 
\begin{equation*}I_\varepsilon(t):=\int_{\R^N}u(t,x)\langle \frac x{R(\varepsilon)}\rangle^{-q}dx.\end{equation*} 
Here, $q$ is a constant satisfying $N<q<N+p\sigma$ and $R(\varepsilon)$ is defined in (\ref{R}).

Multiplying the equation by $\langle \frac x{R(\varepsilon)}\rangle^{-q}$, it follows from H\"{o}lder's inequality and Young's inequality  that
\begin{equation}
\begin{aligned}
& \frac {d^2}{dt^2}I_\varepsilon(t)+b(t)\frac{d}{dt}I_\varepsilon(t)\\
&=\int_{\R^N}\left(\partial_t^2+b(t)\partial_t\right)u(t,x)\langle \frac x{R(\varepsilon)}\rangle^{-q}dx\\
&=-\int_{\R^N}(-\Delta)^{\sigma/ 2}u(t,x)\langle \frac x{R(\varepsilon)}\rangle^{-q}dx+\int_{\R^N}f(u)\langle \frac x{R(\varepsilon)}\rangle^{-q}dx\\
&\geq -\int_{\R^N}(-\Delta)^{\sigma/ 2}u(t,x)\langle \frac x{R(\varepsilon)}\rangle^{-q}dx+C_1\int_{\R^N}|u|^p\langle \frac x{R(\varepsilon)}\rangle^{-q}dx\\
&=-\int_{\R^N} u(t,x)(-\Delta)^{\sigma/ 2}\langle \frac x{R(\varepsilon)}\rangle^{-q}dx+C_1\|u(t)\langle \frac x{R(\varepsilon)}\rangle^{-\frac qp}\|_{L^p}^p\\
&\geq -\|(-\Delta)^{\sigma/ 2}\langle \frac x{R(\varepsilon)}\rangle^{-q}\langle \frac x{R(\varepsilon)}\rangle^{\frac qp}\|_{L^{p'}}\|u(t)\langle \frac x{R(\varepsilon)}\rangle^{\frac q p}\|_{L^p}+C_1\|u(t)\langle \frac x{R(\varepsilon)}\rangle^{-\frac qp}\|_{L^p}^p\\
&\geq -\left(\frac 2{C_1}\right)^{\frac {p'}p}(p')^{-1}p^{-\frac {p'}p}\|(-\Delta)^{\sigma/ 2}\langle \frac x{R(\varepsilon)}\rangle^{-q}\langle \frac x{R(\varepsilon)}\rangle^{\frac qp}\|_{L^{p'}}^{p'}-\frac {C_1}2 \|u(t)\langle \frac x{R(\varepsilon)}\rangle^{-\frac qp}\|_{L^p}^p\\
&\quad \quad +C_1\|u(t)\langle \frac x{R(\varepsilon)}\rangle^{-\frac qp}\|_{L^p}^p\\
&= -\left(\frac 2{C_1}\right)^{p'-1}(p')^{-1}p^{1-p'}\|(-\Delta)^{\sigma/ 2}\langle \frac x{R(\varepsilon)}\rangle^{-q}\langle \frac x{R(\varepsilon)}\rangle^{\frac qp}\|_{L^{p'}}^{p'}+\frac {C_1}2 \|u(t)\langle \frac x{R(\varepsilon)}\rangle^{-\frac qp}\|_{L^p}^p\\
\end{aligned}
\end{equation}

Applying H\"{o}lder's inequality, we have
\begin{equation}
    \begin{aligned}
        I_\varepsilon(t)&=\int_{\R^N}u(t,x)\langle \frac x{R(\varepsilon)}\rangle^{-q}dx\\
        &\leq \left(\int_{\R^N}|u|^p\langle \frac x{R(\varepsilon)}\rangle^{-q}dx\right)^{\frac 1p}\left(\int_{\R^N}\langle \frac x{R(\varepsilon)}\rangle^{-q}dx\right)^{\frac 1{p'}},
    \end{aligned}
\end{equation}
which implies
\begin{equation}
    \|u(t)\langle \frac x{R(\varepsilon)}\rangle^{-\frac qp}\|_{L^p}^p\leq \left(\int_{\R^N}\langle \frac x{R(\varepsilon)}\rangle^{-q}dx\right)^{-\frac p{p'}}|I_\varepsilon(t)|^p.
\end{equation}

Thus, 

\begin{equation}
\begin{aligned}
&\frac {d^2}{dt^2}I_\varepsilon(t)+b(t)\frac{d}{dt}I_\varepsilon(t)\\
&\geq -\left(\frac 2{C_1}\right)^{p'-1}(p')^{-1}p^{1-p'}\|(-\Delta)^{\sigma/ 2}\langle \frac x{R(\varepsilon)}\rangle^{-q}\langle \frac x{R(\varepsilon)}\rangle^{\frac qp}\|_{L^{p'}}^{p'}+\frac {C_1}2 \|\langle \frac x{R(\varepsilon)}\rangle^{-q}\|_{L^1}^{1-p}|I_\varepsilon(t)|^p\\
&=\frac {C_1}2 \|\langle \frac x{R(\varepsilon)}\rangle^{-q}\|_{L^1}^{1-p}\left(|I_\varepsilon(t)|^p-A_\varepsilon(N,p,\sigma,q)^p\right)\\
&=\frac {C_1}2 \|\langle \frac x{R(\varepsilon)}\rangle^{-q}\|_{L^1}^{1-p}\left(I_\varepsilon(t)-A_\varepsilon(N,p,\sigma,q)\right)^p.
\end{aligned}
\end{equation}
Here, $A_\varepsilon(N,p,\sigma,q):=\left(\frac 2{C_1}\right)^{p'-1}(p')^{-1/p}p^{\frac {1-p'}p}\|(-\Delta)^{\sigma/ 2}\langle \frac x{R(\varepsilon)}\rangle^{-q}\langle \frac x{R(\varepsilon)}\rangle^{\frac qp}\|_{L^{p'}}^{\frac {p'}p}\|\langle \frac x{R(\varepsilon)}\rangle^{-q}\|_{L^1}^{\frac 1{p'}}$. Note that the last inequality holds under the condition that 
\begin{equation}\label{CI}
I_\varepsilon(t)-A_\varepsilon(N,p,\sigma,q)>0.
\end{equation}

Introducing $J(t)=I_\varepsilon(t)-A_\varepsilon(N,p,\sigma,q)$, we obtain
\begin{equation}
\left\{
\begin{aligned}
&\frac {d^2}{dt^2}J(t)+b(t)\frac{d}{dt}J(t)\geq \frac {C_1}2 \|\langle \frac x{R(\varepsilon)}\rangle^{-q}\|_{L^1}^{1-p}J(t)^p,\quad \mbox{for}\quad t\in[0,t_0),\\
&J(0) = I_\varepsilon(0)-A_\varepsilon(N,p,\sigma,q)\\
&J'(0) = A_1 J(0)
\end{aligned}
\right.
\end{equation}
where $A_1(\varepsilon):=\frac {I_\varepsilon'(0)}{I_\varepsilon(0)-A_\varepsilon(N,p,\sigma,q)}$.

Direct calculation implies
\begin{equation}\label{Re1}
    \|\langle \frac x{R(\varepsilon)}\rangle^{-q}\|_{L^1}= R(\varepsilon)^{N}\|\langle x\rangle^{-q}\|_{L^1}.
\end{equation}
It follows from Lemma 3.2 in \cite{QD} that
\begin{equation}
    (-\Delta)^{\sigma/ 2}\phi_\varepsilon(x)=\varepsilon^{-\sigma} \left((-\Delta)^{\sigma/ 2}\phi\right)\left(\frac x{\varepsilon}\right),
\end{equation}
with $\phi=\langle x\rangle^{-q}$ and $\phi_\varepsilon(x):=\phi(\frac{x}{\varepsilon})$.
Therefore, we get
\begin{equation}\label{Re2}
    \|(-\Delta)^{\sigma/ 2}\langle \frac x{R(\varepsilon)}\rangle^{-q}\langle \frac x{R(\varepsilon)}\rangle^{\frac qp}\|_{L^{p'}}= R(\varepsilon)^{-\sigma+\frac N{p'}}\|(-\Delta)^{\sigma/ 2}\langle x\rangle^{-q}\langle x\rangle^{\frac qp}\|_{L^{p'}}
\end{equation}
Comparing (\ref{Re1}), (\ref{Re2}), (\ref{A}), (\ref{R}) and the definition of $A_\varepsilon(N,p,\sigma,q)$, we obtain
\begin{equation}\label{AeA}
    A_\varepsilon(N,p,\sigma,q)=A(N,p,\sigma,q)R(\varepsilon)^{-\frac{\sigma p'}p+N}=\frac \varepsilon 4 I_0.
\end{equation}
From the assumption (\ref{BLC1}) and (\ref{AeA}), we obtain
\begin{equation}\label{3.11}
    I_\varepsilon(0)-A_\varepsilon(N,p,\sigma,q)=\varepsilon \int_{\R^N}\langle \frac x {R(\varepsilon_0)}\rangle^{-q}a_0(x)dx-A_\varepsilon(N,p,\sigma,q)\geq \frac \varepsilon 4 I_0.
\end{equation}
And
\[
I_\varepsilon'(0)=\varepsilon\int_{\R^N}\langle \frac x {R(\varepsilon_0)}\rangle^{-q}a_1(x)dx < \varepsilon\int_{\R^N} a_1(x)dx,
\]
which implies that
\begin{equation}
A_1(\varepsilon)\le \frac {4 I_1}{I_0}.
\end{equation}
On the other hand, by the assumption (\ref{BLC3}), we have
\begin{equation}
\begin{aligned}\label{I-A}
I_\varepsilon(0)-A_\varepsilon(N,p,\sigma,q)&=\varepsilon \int_{\R^N}\langle \frac x {R(\varepsilon)}\rangle^{-q}a_0(x)dx-\frac \varepsilon 4 I_0\\
&\leq \varepsilon I_0\\
&\leq \left(\frac {C_1} 2\right)^{-\frac 1{p-1}}\|\langle x\rangle^{-q}\|_{L^1}R(\varepsilon)^N\\
&=\left(\frac {C_1} 2\right)^{-\frac 1{p-1}}\|\langle \frac x {R(\varepsilon)}\rangle^{-q}\|_{L^1}\\
\end{aligned}
\end{equation}
since $\int_{\R^N}\langle \frac x {R(\varepsilon)}\rangle^{-q}a_0(x)dx\rightarrow \int_{\R^N}a_0(x)dx=I_0$ as $\varepsilon\rightarrow 0$. Moreover, by the assumption (\ref{BLC2}) and (\ref{I-A}),
\begin{equation}\label{3.13}
\begin{aligned}
     A_1(\varepsilon)&=\frac {I_\varepsilon'(0)}{I_\varepsilon(0)-A_\varepsilon(N,p,\sigma,q)}\\
     &= \frac {\varepsilon\int_{\R^N}\langle \frac x {R(\varepsilon_0)}\rangle^{-q}a_1(x)dx}{I_\varepsilon(0)-A_\varepsilon(N,p,\sigma,q)}\\
     &\geq \frac {\frac \varepsilon 2 I_1}{\varepsilon I_0}=\frac {I_1}{2 I_0}
\end{aligned}
\end{equation}

Now, from (\ref{3.11})-(\ref{3.13}), the conditions in Proposition 2.3 in \cite{FIW2019} are fulfilled. Thus, applying this proposition, we get
\begin{equation}\label{estimateJ}
J(t)\geq J(0)(1-\mu(p, b,\beta,A_1)\Tilde{J}(0)^{p-1}B(t))^{-\frac 2{p-1}},
\end{equation}
and the upper bound of the lifespan
\begin{equation}\label{estimateT}
    T_0\leq B^{-1}(\mu(p,b,\beta, A_1)^{-1}\Tilde{J}(0)^{1-p}),
\end{equation}
where $\Tilde{J}(0)=\left(\frac {C_2}2\right)^{\frac 1{p-1}} \|\langle \frac x{R(\varepsilon)}\rangle^{-q}\|_{L^1}^{-1}J(0)$, and
\begin{equation*}
\begin{aligned}
\mu(p,b,\beta, A_1)&=\min\left\{1,\frac{p-1}2 b(0)A_1\right., \\
&\quad \quad \quad \left.\left[\frac {2(p+1)}{(p-1)^2}b_1^{-2}(2^{\frac 1{1+\beta}}(1+B_4))^{\max(0,2\beta)}+\frac {2(b_1^{-1}b_3+1)}{p-1}\right]^{-1}\right\}.
\end{aligned}
\end{equation*}
Here, $B(t):=\int_0^tb(s)^{-1}ds$, and $B^{-1}$ is the inverse function of $B(t)$, satisfying 
\begin{equation}
\begin{cases}
    B_1 (1+t)^{1/(1+\beta)}\leq B^{-1}(t) \leq B_2 (1+t)^{1/(1+\beta)},& \quad \mbox{for}\quad -1<\beta\leq 1,\\
    \exp{(B_1 (1+t))}\leq B^{-1}(t) \leq \exp{(B_2 (1+t))},& \quad \mbox{for}\quad \beta=-1,
\end{cases}
\end{equation}
with some constants $B_1,B_2>0$.

Now, let's prove the condition (\ref{CI}) holds, that is, $J(t)>0$ for all $t\in [0,T_0)$. Let $T^*>0$ satisfy $J(t)>0$ for $t\in [0,T^*)$ and $J(T^*)=0$. Then, in the same manner as above, we have (\ref{estimateJ}) for $t\in [0,T_0)$. However, the right hand side of (\ref{estimateJ}) is still positive for $t=T^*$, which contradicts $J(T^*)=0$. Thus, we could conclude that $J(t)>0$ for all $t\in [0,T_0)$, which implies that (\ref{estimateJ}) and (\ref{estimateT}) hold for $t\in [0,T_0)$.

Note that
\begin{equation}\label{estimateJt}
\begin{aligned}
\Tilde{J}(0)^{p-1}&=\frac {C_1}2 \|\langle \frac x{R(\varepsilon)}\rangle^{-q}\|_{L^1}^{1-p}J(0)^{p-1}\\
&\geq \frac {C_1}2 \|\langle x\rangle^{-q}\|_{L^1}^{1-p}R(\varepsilon)^{N(1-p)}\left(\frac \varepsilon 4 I_0\right)^{p-1}\\
&=\frac {C_1}2 \|\langle x\rangle^{-q}\|_{L^1}^{1-p}A(N,p,\sigma,q)^{\frac{(p-1)^2}{p-p_c}}\left(\frac \varepsilon 4 I_0\right)^{-\frac{(p-1)^2}{p-p_c}+p-1}\\
&=\frac {C_1}2 \|\langle x\rangle^{-q}\|_{L^1}^{1-p}A(N,p,\sigma,q)^{\frac{(p-1)^2}{p-p_c}}\left(\frac \varepsilon 4 I_0\right)^{\frac{1}{\frac 1{p-1}-\frac N\sigma}},
\end{aligned}
\end{equation}
since
\begin{equation*}
-\frac{(p-1)^2}{p-p_c}+p-1=\frac{1}{\frac 1{p-1}-\frac N\sigma} 
\end{equation*}
by direct calculation. Thus, substituting the estimate (\ref{estimateJt}) into (\ref{estimateT}), we obtain
\begin{equation}
T_0\leq 
\begin{cases}
C \varepsilon^{-\frac 1{\left(\frac 1{p-1}-\frac N\sigma\right)(1+\beta)}}, &\quad \mbox{for}\quad -1<\beta\leq 1,\\
\exp{\left(C\varepsilon^{-\frac{1}{\frac 1{p-1}-\frac N\sigma}}\right)}, &\quad \mbox{for}\quad \beta= -1.
\end{cases}
\end{equation} 

\section{Proof of blow-up and upper estimate of lifespan in the critical case}

To rewrite the equation into a divergence form, we introduce a function \(g:[0,\infty)\rightarrow \R_{>0}\) as the unique solution to the initial value problem
of the ordinary differential equation
\begin{equation}\label{gtODE}
 \left\{ \begin{array}{l}
  \displaystyle{-g'(t)+b(t)g(t)=1, \ \ t\in (0,\infty), }\\
  \displaystyle{g(0)=B_0,}
 \end{array}\right.
\end{equation}
where
\(B_0:=\int_0^{\infty} \exp{\left(-\int_0^t b(s)\,ds\right)}\,dt<\infty.\)
The solution \(g\) is explicitly given by
\begin{equation}\label{gt}
 g(t) = \exp{\left(\int_0^t b(s)\,ds\right)} \left( B_0 - \int_0^t \exp{\left(-\int_0^{\tau} b(s)\,ds\right)}\,d\tau \right).
\end{equation}
Here we note that \(g'(0)=b(0)g(0)-1\) and
\begin{equation}\label{2.5}
 \lim_{t \to \infty}b(t)g(t) = 1, \ \ C^{-1}b(t)^{-1} \le g(t) \le Cb(t)^{-1} \  (0 \le t <\infty)
\end{equation}

With the property of the function \(g\), we can get the following equation of divergence form.
\begin{lemma}
Let \(u\) be a strong solution. Then \(u\) satisfies
\begin{equation}\label{DF}
 (g(t)u)_{tt} +(-\Delta)^{\sigma/2} (g(t)u) -(g'(t)u)_t + u_t = g(t)f(u)
\end{equation}
in the sense of \(C([0,T);L^2(\R^N))\)
\end{lemma}

Set \(F(x,t):=g(t)f(u)-(g(t)u)_{tt}+(g'(t)u)_t\). Thus, we could write (\ref{DF}) as
\begin{equation*}u_t+(-\Delta)^{\sigma/2} (g(t)u) =F(x,t).\end{equation*}
Taking the Fourier transform, we get
\begin{equation*}\partial_t \hat{u}+g(t)|\xi|^{\sigma} \hat{u} =\hat{F}(t).\end{equation*}
Solving this ordinary differential equation, we have
\begin{equation}\label{fouriersol}
\hat{u}(t) = e^{-G(t)|\xi|^{\sigma}} \hat{u}(0)+\int^t_0 e^{-(G(t)-G(s))|\xi|^{\sigma}} F(s)ds,
\end{equation}
where \(G(t):=\int_0^t g(s)ds\).

Denote
\begin{equation*}\varphi(x):=(2\pi)^{-N/2}\mathcal{F}^{-1}[e^{-|\xi|^\sigma}]=\frac 1{(2\pi)^{N}}\int_{\R^N}e^{-|\xi|^\sigma}e^{i\xi\cdot x}d\xi,\end{equation*}
and
\begin{equation*}\Phi(t)=\Phi(x,t):=\frac 1{t^{N/\sigma}}\varphi \left(\frac x{t^{1/\sigma}}\right)=(2\pi)^{-N/2}\mathcal{F}^{-1}[e^{-t|\xi|^\sigma}].\end{equation*}
Then, taking the inverse Fourier transform of (\ref{fouriersol}), we obtain
\begin{lemma}
Let \(u\) be a strong solution. Then
\begin{equation}\label{integralsol}
u(t)=\Phi(G(t))\ast u_0 +\int_0^t \Phi(G(t)-G(s))\ast F(s) ds.
\end{equation}

\end{lemma}

To prove blow-up in the critical case, we require the following estimates of \(\varphi(x)\).
\begin{lemma}\label{buc1}
For \(\sigma\in (0,2)\), the following estimates 
\begin{equation}\label{buc1-1}
\int_{\R^N}\varphi(x)dx = 1,
\end{equation}
\begin{equation}\label{buc1-2}
\int_{\R^N}|\nabla \varphi(x)\cdot x|dx \leq C,
\end{equation}
and
\begin{equation}\label{buc1-3}
\int_{\R^N}|x|^2 |\Delta \varphi(x)|dx \leq C,
\end{equation}
hold for some constant \(C\).
\end{lemma}

{\bf Proof.} Since
\[
\widehat{\varphi}(0) = (2\pi)^{-N/2}\int_{\mathbb R^N}\varphi(x)\,dx,
\]
we have
\[
\int_{\mathbb R^N}\varphi(x)\,dx=1.
\] 

It follows from Theorem 2.1 of \cite{KL} that 
\begin{equation}
|\nabla^i \varphi(x)|\leq \frac{C}{|x|^{N+\sigma+i}}
\end{equation}
for \(|x|\geq 1\) (\(i=1,2\)). Thus,
\begin{equation}
\begin{aligned}
\int_{\R^N} \nabla \varphi(x)dx &= \int_{|x|<1} \nabla \varphi(x)dx+\int_{|x|\geq 1} \nabla \varphi(x)dx\\
&\leq C+C\int_{|x|\geq 1} |x|^{-(N+\sigma+1)}dx\\
&\leq C
\end{aligned}
\end{equation}
The inequality (\ref{buc1-3}) could be obtained similarly.

{\bf Proof of Theorem 1.3.} Multiplying the equation (\ref{integralsol}) by \((G(t)+1)^{N/\sigma}\Phi(x,G(t)+1)\) and integrating it on \(\R^N\), we get
\begin{equation}
\begin{aligned}
&\int_{\R^N}\varphi\left(\frac x{(G(t)+1)^{1/\sigma}}\right)u(x,t)dx\\
&=(G(t)+1)^{N/\sigma}\int_{\R^N}\Phi(x,G(t)+1)\Phi(G(t))\ast u_0 dx\\
&+(G(t)+1)^{N/\sigma}\int_{\R^N}\Phi(x,G(t)+1)\int_0^t \Phi(G(t)-G(s))\ast g(s)f(u) ds dx\\
&-(G(t)+1)^{N/\sigma}\int_{\R^N}\Phi(x,G(t)+1)\int_0^t \Phi(G(t)-G(s))\ast (g(s)u(s))_{ss} ds dx\\
&+(G(t)+1)^{N/\sigma}\int_{\R^N}\Phi(x,G(t)+1)\int_0^t \Phi(G(t)-G(s))\ast (g'(s)u(s))_s ds dx\\
\end{aligned}
\end{equation}
It follows from the property \(\Phi(s)\ast \Phi(t) =\Phi(s+t)\) for \(s,t>0\), and Fubini-Tonelli theorem that
\begin{equation}
\int_{\R^N}\Phi(x,G(t)+1)\Phi(G(t))\ast u_0 dx
=\int_{\R^N}\Phi(2G(t)+1) u_0(x)dx,\\
\end{equation}
and
\begin{equation}
\begin{aligned}
&\int_{\R^N}\Phi(x,G(t)+1)\int_0^t \Phi(G(t)-G(s))\ast (g(s)u(s))_{ss} ds dx\\
&-\int_{\R^N}\Phi(x,G(t)+1)\int_0^t \Phi(G(t)-G(s))\ast (g'(s)u(s))_s ds dx\\
&=\int_{\R^N}\Phi(x,G(t)+1)\int_0^t \Phi(G(t)-G(s))\ast (g(s)u_s(s))_s ds dx\\
&=\int_0^t\int_{\R^N}\Phi(x,2G(t)-G(s)+1)\partial_s (g(s)u_s(s)) ds dx\\
&=\int_{\R^N} \Phi(x,G(t)+1) g(t)u_t(x,t)dx - \int_{\R^N} \Phi(x,2G(t)+1) g(0)u_t(x,0)dx\\
&\ \ - \int_0^t\int_{\R^N} \left[ \partial_s \Phi(x,2G(t)-G(s)+1) \right] g(s)u_s(s) ds dx\\
&=\int_{\R^N} \Phi(x,G(t)+1) g(t)u_t(x,t)dx - \int_{\R^N} \Phi(x,2G(t)+1) g(0)u_t(x,0)dx\\
&\ \ +\int_0^t\int_{\R^N}g(s)^2\Phi_2'(x,2G(t)-G(s)+1)u_s(s) ds dx\\
\end{aligned}
\end{equation}
where \(\Phi_2'(x,t):=\frac{\partial}{\partial t}\Phi(x,t)\).

Thus,
\begin{equation}\label{11}
\begin{aligned}
&\int_{\R^N}\varphi\left(\frac x{(G(t)+1)^{1/\sigma}}\right)u(x,t)dx\\
&=(G(t)+1)^{N/\sigma}\int_{\R^N}\Phi(2G(t)+1) u_0(x)dx\\
&+(G(t)+1)^{N/\sigma}\int_0^t\int_{\R^N}\Phi(2G(t)-G(s)+1)g(s)f(u(s)) dxds\\
&-(G(t)+1)^{N/\sigma}\int_{\R^N} \Phi(x,G(t)+1) g(t)u_t(x,t)dx\\
&+(G(t)+1)^{N/\sigma} \int_{\R^N} \Phi(x,2G(t)+1) g(0)u_t(x,0)dx\\
&-(G(t)+1)^{N/\sigma}\int_0^t\int_{\R^N}g(s)^2\Phi_2'(x,2G(t)-G(s)+1)u_s(s) ds dx\\
\end{aligned}
\end{equation}
Let
\begin{equation*}
\begin{aligned}
\mathcal{A}(t)&:=\int_{\R^N}\varphi\left(\frac x{(G(t)+1)^{1/\sigma}}\right)u(x,t)dx\\
\mathcal{B}(t)&:=(G(t)+1)^{N/\sigma}\int_{\R^N} \Phi(x,G(t)+1) g(t)u_t(x,t)dx\\
\mathcal{C}(t)&:=(G(t)+1)^{N/\sigma}\int_{\R^N}\Phi(2G(t)+1) (u_0+g(0)u_1)(x)dx\\
\mathcal{D}(t)&:=(G(t)+1)^{N/\sigma}\int_0^t\int_{\R^N}\Phi(2G(t)-G(s)+1)g(s)f(u(s)) dxds\\
\mathcal{E}(t)&:=-(G(t)+1)^{N/\sigma}\int_0^t\int_{\R^N}g(s)^2\Phi_2'(x,2G(t)-G(s)+1)u_s(s) ds dx\\
\end{aligned}
\end{equation*}
Then, the equation (\ref{11}) becomes \(\mathcal{A}(t)+\mathcal{B}(t)=\mathcal{C}(t)+\mathcal{D}(t)+\mathcal{E}(t)\).

We estimate \(\mathcal{A}(t)\) first. It follows from H\"{o}lder's inequality, we have
\begin{equation}\label{estimateA}
\begin{aligned}
\mathcal{A}(t)&=\int_{\R^N}\varphi\left(\frac x{(G(t)+1)^{1/\sigma}}\right)u(x,t)dx\\
&\leq \left(\int_{\R^N} \varphi(\frac x{(G(t)+1)^{1/\sigma}})dx\right)^{1/p'}\left(\int_{\R^N} \varphi(\frac x{(G(t)+1)^{1/\sigma}})|u(t,x)|^pdx\right)^{1/p}\\
&\leq C (G(t)+1)^{N/\sigma p'}\left(\int_{\R^N} \varphi(\frac x{(G(t)+1)^{1/\sigma}})|u(t,x)|^pdx\right)^{1/p}=:C\mathcal{R}(t),\\
\end{aligned}
\end{equation}
where \(p'\) is the H\"{o}lder conjugate of \(p\) with \(p'=\frac p{p-1}\), since \(\int_{\R^N} \varphi(x) dx \leq C\) by Lemma \ref{buc1}.

And
\begin{equation*}
\begin{aligned}
\mathcal{B}(t)&=(G(t)+1)^{N/\sigma}\int_{\R^N} \Phi(x,G(t)+1) g(t)u_t(x,t)dx\\
&=g(t)\mathcal{A}'(t)+\frac {g(t)^2}{\sigma (G(t)+1)}\int_{\R^N}\nabla \varphi(\frac x{(G(t)+1)^{1/\sigma}})\cdot \frac x{(G(t)+1)^{1/\sigma}} u(x,t)dx.
\end{aligned}
\end{equation*}
Notice that 
\begin{equation}
\varphi(x)\sim |x|^{-N-\sigma}\quad \mbox{and}\quad |\nabla \varphi|\le C|x|^{-N-\sigma-1}
\end{equation}
as $x\rightarrow \infty$, from Theorem 2.1 in \cite{KL}. Therefore, by H\"{o}lder inequality,
\begin{equation}\label{B2}
\begin{aligned}
&\left|\frac {g(t)^2}{\sigma (G(t)+1)}\int_{\R^N}\nabla \varphi(\frac x{(G(t)+1)^{1/\sigma}})\cdot \frac x{(G(t)+1)^{1/\sigma}} u(x,t)dx\right|\\
&\leq \frac {g(t)^2}{\sigma (G(t)+1)} \left(\int_{\R^N}\left|\nabla \varphi(\frac x{(G(t)+1)^{1/\sigma}})\cdot \frac x{(G(t)+1)^{1/\sigma}}\right|^{p'}\left|\varphi(\frac x{(G(t)+1)^{1/\sigma}})\right|^{-\frac {p'}p}dx\right)^{1/p'}\\
&\quad\quad \times\left(\int_{\R^N} \varphi(\frac x{(G(t)+1)^{1/\sigma}})|u(t,x)|^pdx\right)^{1/p}\\
&\leq C \frac {g(t)^2}{G(t)+1} (G(t)+1)^{N/\sigma p'}\left(\int_{\R^N} \varphi(\frac x{(G(t)+1)^{1/\sigma}})|u(t,x)|^pdx\right)^{1/p}=C \frac {g(t)^2}{G(t)+1} \mathcal{R}(t),
\end{aligned}
\end{equation}
which implies
\begin{equation}
\label{estimateB}
\mathcal{B}(t)\leq g(t)\mathcal{A}'(t)+C \frac {g(t)^2}{G(t)+1} \mathcal{R}(t).
\end{equation}

We estimate \(\mathcal{D}(t)\) as follows.
\begin{equation*}
\begin{aligned}
\mathcal{D}(t)&\geq C_1(G(t)+1)^{N/\sigma}\int_0^t\int_{\R^N}\Phi(2G(t)-G(s)+1)g(s)|u(x,s)|^p dxds\\
&=C_1\int_0^t\left(\frac {G(t)+1}{2G(t)-G(s)+1}\right)^{N/\sigma}g(s)\int_{\R^N}\varphi(\frac x{(2G(t)-G(s)+1)^{1/\sigma}})|u(x,s)|^pdxds\\
&\geq 2^{-N/\sigma}C_1\int_0^t g(s)\int_{\R^N}\varphi(\frac x{(2G(t)-G(s)+1)^{1/\sigma}})|u(x,s)|^pdxds\\
&= 2^{-N/\sigma}C_1\int_0^t g(s)(G(s)+1)^{-\frac {N}{\sigma}(p-1)}\mathcal{R}(s)^pds\\
\end{aligned}
\end{equation*}

Then, for \(\mathcal{E}(t)\), by the integration by parts, we have
\begin{equation*}
\begin{aligned}
\mathcal{E}(t)&=-(G(t)+1)^{N/\sigma}\int_0^t\int_{\R^N}g(s)^2\Phi_2'(x,2G(t)-G(s)+1)u_s(s) ds dx\\
&=-(G(t)+1)^{N/\sigma}\int_{\R^N}g(t)^2\Phi_2'(x,G(t)+1)u(t,x) dx\\
&\ +(G(t)+1)^{N/\sigma}\int_{\R^N}g(0)^2\Phi_2'(x,2G(t)+1)u_0(x) dx\\
&\ +(G(t)+1)^{N/\sigma}\int_0^t\int_{\R^N}\left[g(s)^2\Phi_2'(x,2G(t)-G(s)+1)\right]_s u(x,s) ds dx\\
&=-(G(t)+1)^{N/\sigma}\int_{\R^N}g(t)^2\Phi_2'(x,G(t)+1)u(t,x) dx\\
&\ +(G(t)+1)^{N/\sigma}\int_{\R^N}g(0)^2\Phi_2'(x,2G(t)+1)u_0(x) dx\\
&\ +2(G(t)+1)^{N/\sigma}\int_0^t\int_{\R^N}g(s)g'(s)\Phi_2'(x,2G(t)-G(s)+1)u(x,s) ds dx\\
&\ -(G(t)+1)^{N/\sigma}\int_0^t\int_{\R^N}g(s)^3\Phi_{22}''(x,2G(t)-G(s)+1) u(x,s) ds dx\\
&=:\mathcal{E}_1(t)+\mathcal{E}_2(t)+\mathcal{E}_3(t)+\mathcal{E}_4(t).
\end{aligned}
\end{equation*}
where \(\Phi_{22}''(x,t):=\frac{\partial^2}{\partial t^2}\Phi(x,t)\).

Since
\begin{equation}\label{DPhi}
\Phi_{2}'(x,t)=\frac {\partial} {\partial t}\Phi(x,t) = -\frac N{\sigma t} \Phi(x,t)-\frac 1{\sigma t^{N/\sigma+1}} \nabla \varphi(\frac x{t^{1/\sigma}})\cdot \frac x{t^{1/\sigma}},
\end{equation}
we get
\begin{equation}\label{E1}
\begin{aligned}
|\mathcal{E}_1(t)|&\leq (G(t)+1)^{N/\sigma} \frac{N g(t)^2}{\sigma (G(t)+1)}\int_{\R^N}|\Phi(x,G(t)+1) u(x,t)|dx\\
&\ +\left|\frac {g(t)^2}{\sigma (G(t)+1)}\int_{\R^N}\nabla \varphi(\frac x{(G(t)+1)^{1/\sigma}})\cdot \frac x{(G(t)+1)^{1/\sigma}} u(x,t)dx\right|\\
&\leq C \frac {g(t)^2}{G(t)+1} \mathcal{R}(t)
\end{aligned}
\end{equation}
by (\ref{estimateA}) and (\ref{B2}), and
\begin{equation}\label{E2}
\begin{aligned}
\mathcal{E}_2(t)&= -\left(\frac{G(t)+1}{2G(t)+1}\right)^{N/\sigma} \frac{N g(0)^2}{\sigma (2G(t)+1)}\int_{\R^N}|\varphi(\frac x{(2G(t)+1)^{1/\sigma}}) u_0(x)|dx\\
&\ +\left(\frac{G(t)+1}{2G(t)+1}\right)^{N/\sigma}\frac {g(0)^2}{2\sigma (G(t)+1)}\int_{\R^N}\nabla \varphi(\frac x{(2G(t)+1)^{1/\sigma}})\cdot \frac x{(2G(t)+1)^{1/\sigma}} u_0(x)dx
\end{aligned}
\end{equation}
hold.

It follows from (\ref{DPhi}) that
\begin{equation}\label{estimateE3}
\begin{aligned}
&|\mathcal{E}_3(t)|\\
&\leq (G(t)+1)^{N/\sigma} \frac{2N}{\sigma (G(t)+1)}\int_{\R^N}\int_0^t g(s)|g'(s)||\Phi(x,2G(t)-G(s)+1) u(x,s)|dsdx\\
&\ +\frac 2\sigma (G(t)+1)^{\frac N\sigma-1}\int_0^t\int_{\R^N}(2G(t)-G(s)+1)^{-\frac {N+1}\sigma}g(s)|g'(s)| |\nabla \varphi(\frac x{(2G(t)-G(s)+1)^{1/\sigma}})\cdot x||u(x,s)|dxds\\
&=:\mathcal{E}_{31}(t)+\mathcal{E}_{32}(t).
\end{aligned}
\end{equation}

By H\"{o}lder inequality and Young's inequality, we get
\begin{equation}\label{E31}
\begin{aligned}
\mathcal{E}_{31}(t)&\leq \frac{2N}{\sigma}(G(t)+1)^{N/(\sigma p')-1}\left(\int_{\R^N}\int_0^t g(s)|g'(s)|^{p'}|\Phi(x,2G(t)-G(s)+1)|dsdx\right)^{1/p'}\\
& \quad \times (G(t)+1)^{N/(\sigma p)}\left(\int_0^t\int_{\R^N}\Phi(x,2G(t)-G(s)+1)g(s)|u(s)|^p dxds\right)^{1/p}\\
&\leq C (G(t)+1)^{N/(\sigma p')-1}\left(\int_0^t g(s)|g'(s)|^{p'}ds\right)^{1/p'} \mathcal{D}(t)^{1/p}\\
&\leq \frac 18 \mathcal{D}(t)+C(G(t)+1)^{N/(\sigma)-p'}\left(\int_0^t g(s)|g'(s)|^{p'}ds\right)
\end{aligned}
\end{equation}
and
\begin{equation}\label{E32}
\begin{aligned}
&\mathcal{E}_{32}(t)\leq \frac{2}{\sigma}(G(t)+1)^{N/\sigma-1}\left(\int_0^t\int_{\R^N}\Phi(x,2G(t)-G(s)+1)g(s)|u(s)|^p dxds\right)^{1/p}\\
&\times \left(\int_0^t\int_{\R^N}(2G(t)-G(s)+1)^{-\frac {N+p'}\sigma}g(s)|g'(s)|^{p'} |\nabla \varphi(\frac x{(2G(t)-G(s)+1)^{1/\sigma}})\cdot x|^{p'}|\varphi(\frac x{(2G(t)-G(s)+1)^{1/\sigma}})|^{-\frac {p'}p}dxds\right)^{1/p'}\\
&\leq C (G(t)+1)^{N/(\sigma p')-1}\mathcal{D}(t)^{1/p}\left(\int_0^t g(s)|g'(s)|^{p'}ds\right)^{1/p'}\\
&\leq \frac 18 \mathcal{D}(t)+C(G(t)+1)^{N/\sigma-p'}\left(\int_0^t g(s)|g'(s)|^{p'}ds\right)
\end{aligned}
\end{equation}
Combining the estimates (\ref{estimateE3})-(\ref{E32}), we obtain
\begin{equation}\label{E3}
|\mathcal{E}_3(t)|\leq \frac 14 \mathcal{D}(t)+C(G(t)+1)^{N/\sigma-p'}\left(\int_0^t g(s)|g'(s)|^{p'}ds\right).
\end{equation}

Direct computation shows
\begin{equation*}
\begin{aligned}
\Phi_{22}''(x,t)&=\frac {N(N+1)}{\sigma^2 t^2}\Phi(x,t)+\frac {2N+\sigma+1}{\sigma^2 t^{2+N/\sigma}}\nabla \varphi(\frac x{t^{1/\sigma}})\cdot \frac x{t^{1/\sigma}}\\
&\ +\frac 1{\sigma t^{2+N/\sigma}} \Delta \varphi(\frac x{t^{1/\sigma}}) \frac{|x|^2}{t^{2/\sigma}},
\end{aligned}
\end{equation*}
which implies
\begin{equation}\label{estimateE4}
\begin{aligned}
&|\mathcal{E}_4(t)|\\
&=(G(t)+1)^{N/\sigma}\int_0^t\int_{\R^N}g(s)^3|\Phi_{22}''(x,2G(t)-G(s)+1) u(x,s)| ds dx\\
&\leq C(G(t)+1)^{N/\sigma-2}\int_0^t\int_{\R^N}g(s)^3|\Phi(x,2G(t)-G(s)+1) u(x,s)| ds dx\\
&\ +C(G(t)+1)^{N/\sigma-2}\int_0^t\int_{\R^N}g(s)^3(2G(t)-G(s)+1)^{-\frac {N+1}\sigma}g(s)|g'(s)| |\nabla \varphi(\frac x{(2G(t)-G(s)+1)^{1/\sigma}})\cdot x||u(x,s)|dxds\\
&\ +C(G(t)+1)^{N/\sigma-2}\int_0^t\int_{\R^N}g(s)^3(2G(t)-G(s)+1)^{-\frac {N+1}\sigma}g(s)|g'(s)| |\Delta \varphi(\frac {x}{(2G(t)-G(s)+1)^{2/\sigma}}) |x|^2 | |u(x,s)|dxds\\
&=:\mathcal{E}_{41}(t)+\mathcal{E}_{42}(t)+\mathcal{E}_{43}(t).
\end{aligned}
\end{equation}
Again, by H\"{o}lder's inequality and Young's inequality, we get
\begin{equation}\label{E41}
\begin{aligned}
\mathcal{E}_{41}(t)&\leq C(G(t)+1)^{N/\sigma-2}\left(\int_0^t\int_{\R^N}g(s)|\Phi(x,2G(t)-G(s)+1)| |u(x,s)|^p dx ds\right)^{1/p}\\
& \quad \times \left(\int_0^t\int_{\R^N}\Phi(x,2G(t)-G(s)+1)g(s)^{2p'+1} dxds\right)^{1/p}\\
&\leq C (G(t)+1)^{N/(\sigma p')-2}\left(\int_0^t g(s)^{2p'+1}ds\right)^{1/p'} \mathcal{D}(t)^{1/p}\\
&\leq \frac 1{12} \mathcal{D}(t)+C(G(t)+1)^{N/\sigma-2p'}\left(\int_0^t g(s)^{2p'+1}ds\right),
\end{aligned}
\end{equation}
\begin{equation}\label{E42}
\begin{aligned}
&\mathcal{E}_{42}(t)\leq C(G(t)+1)^{N/\sigma-2}\left(\int_0^t\int_{\R^N}g(s)|\Phi(x,2G(t)-G(s)+1)| |u(x,s)|^p dx ds\right)^{1/p}\\
&\ \times \left(\int_0^t\int_{\R^N}(2G(t)-G(s)+1)^{-\frac {N+p'}\sigma}g(s)^{2p'+1} |\nabla \varphi(\frac x{(2G(t)-G(s)+1)^{1/\sigma}})\cdot x|^{p'}|\varphi(\frac x{(2G(t)-G(s)+1)^{1/\sigma}})|^{-\frac {p'}p}dxds\right)^{1/p'}\\
&\leq C (G(t)+1)^{N/(\sigma p')-2}\left(\int_0^t g(s)^{2p'+1}ds\right)^{1/p'} \mathcal{D}(t)^{1/p}\\
&\leq \frac 1{12} \mathcal{D}(t)+C(G(t)+1)^{N/\sigma-2p'}\left(\int_0^t g(s)^{2p'+1}ds\right),
\end{aligned}
\end{equation}
and
\begin{equation}\label{E43}
\begin{aligned}
&\mathcal{E}_{43}(t)\leq C(G(t)+1)^{N/\sigma-2}\left(\int_0^t\int_{\R^N}g(s)|\Phi(x,2G(t)-G(s)+1)| |u(x,s)|^p dx ds\right)^{1/p}\\
&\ \times \left(\int_0^t\int_{\R^N}(2G(t)-G(s)+1)^{-\frac {N+2p'}\sigma}g(s)^{2p'+1} |\Delta \varphi(\frac x{(2G(t)-G(s)+1)^{1/\sigma}})|x|^2|^{p'}|\varphi(\frac x{(2G(t)-G(s)+1)^{1/\sigma}})|^{-\frac {p'}p}dxds\right)^{1/p'}\\
&\leq C (G(t)+1)^{N/(\sigma p')-2}\left(\int_0^t g(s)^{2p'+1}ds\right)^{1/p'} \mathcal{D}(t)^{1/p}\\
&\leq \frac 1{12} \mathcal{D}(t)+C(G(t)+1)^{N/\sigma-2p'}\left(\int_0^t g(s)^{2p'+1}ds\right),
\end{aligned}
\end{equation}
since
\begin{equation*}
\int_{\R^N}|\Delta \varphi(x)|x|^2|^{p'}|\varphi(x)|^{-\frac {p'}p}dx\leq C.
\end{equation*}
Thus, combining the estimates (\ref{estimateE4})-(\ref{E43}), we get
\begin{equation}\label{E4}
|\mathcal{E}_4(t)|\leq \frac 1{4} \mathcal{D}(t)+C(G(t)+1)^{N/\sigma-2p'}\int_0^t g(s)^{2p'+1}ds.
\end{equation}

Therefore, from the estimates of \(\mathcal{B}(t), \mathcal{E}(t)\) and \(\mathcal{D}(t)\), we have
\begin{equation}\label{bu31}
\mathcal{A}(t)+g(t)\mathcal{A}'(t)+C \frac {g(t)^2}{G(t)+1}\mathcal{R}(t)\geq C\int_0^t\frac {g(s)\mathcal{R}(s)^p}{(G(s)+1)^{\frac N\sigma (p-1)}}ds + \varepsilon \mathcal{H}(t)-\mathcal{I}(t),
\end{equation}
where \(\mathcal{H}(t) = \varepsilon^{-1}(\mathcal{C}(t)-\mathcal{E}_2(t))\) and
\begin{equation*}
\mathcal{I}(t)=C(G(t)+1)^{N/\sigma-p'}\int_0^t g(s)|g'(s)|^{p'}ds+C(G(t)+1)^{N/\sigma-2p'}\int_0^t g(s)^{2p'+1}ds.
\end{equation*}

Define \(\Gamma(t):=\int_0^t \frac 1{g(s)}ds\). Multiplying both sides of the equation (\ref{bu31}) by \(\left(e^{\Gamma(t)}\right)'\), we get
\begin{equation}
\left(e^{\Gamma(t)}\mathcal{A}(t)\right)'+C \frac {g(t)e^{\Gamma(t)}}{G(t)+1}\mathcal{R}(t)\geq C\left(e^{\Gamma(t)}\right)'\int_0^t\frac {g(s)\mathcal{R}(s)^p}{(G(s)+1)^{\frac N\sigma (p-1)}}ds + \left(e^{\Gamma(t)}\right)'(\varepsilon \mathcal{H}(t)-\mathcal{I}(t)),
\end{equation}
since \(\left(e^{\Gamma(t)}\right)'=\frac {e^{\Gamma(t)}}{g(t)}\).

Integrating on \([0,T]\), 
\begin{equation}\label{bu33}
\begin{aligned}
&e^{\Gamma(t)}\mathcal{A}(t)-\mathcal{A}(0)+C \int_0^t\frac {g(\tau)e^{\Gamma(\tau)}}{G(\tau)+1}\mathcal{R}(\tau)d\tau\\
&\geq C\int_0^t \left(e^{\Gamma(\tau)}\right)'\int_0^\tau\frac {g(s)\mathcal{R}(s)^p}{(G(s)+1)^{\frac N\sigma (p-1)}}dsd\tau + \int_0^t\left(e^{\Gamma(\tau)}\right)'(\varepsilon \mathcal{H}(\tau)-\mathcal{I}(\tau))d\tau\\
&=C\int_0^t \frac {(e^{\Gamma(t)}-e^{\Gamma(s)})g(s)\mathcal{R}(s)^p}{(G(s)+1)^{\frac N\sigma (p-1)}}ds + \int_0^t\left(e^{\Gamma(\tau)}\right)'(\varepsilon \mathcal{H}(\tau)-\mathcal{I}(\tau))d\tau
\end{aligned}
\end{equation}

Combining (\ref{bu33}) with (\ref{estimateA}), we have
\begin{equation}\label{bu34}
\begin{aligned}
&e^{\Gamma(t)}\mathcal{R}(t)+C \int_0^t\frac {g(\tau)e^{\Gamma(\tau)}}{G(\tau)+1}\mathcal{R}(\tau)d\tau\\
&\geq C\int_0^t \frac {(e^{\Gamma(t)}-e^{\Gamma(s)})g(s)\mathcal{R}(s)^p}{(G(s)+1)^{\frac N\sigma (p-1)}}ds + \int_0^t\left(e^{\Gamma(\tau)}\right)'(\varepsilon \mathcal{H}(\tau)-\mathcal{I}(\tau))d\tau+\mathcal{A}(0).
\end{aligned}
\end{equation}

Multiplying \(g(t)^{-1}\) and integrating on \([0,t]\), we obtain
\begin{equation}\label{bu35}
\begin{aligned}
&\int_0^t \frac 1{g(s)} e^{\Gamma(s)}\mathcal{R}(s)ds+C \int_0^t \frac 1{g(s)}\int_0^s\frac {g(\tau)e^{\Gamma(\tau)}}{G(\tau)+1}\mathcal{R}(\tau)d\tau ds\\
&\geq C\int_0^t \frac 1{g(s)}\int_0^s \frac {(e^{\Gamma(s)}-e^{\Gamma(\tau)})g(\tau)\mathcal{R}(\tau)^p}{(G(\tau)+1)^{\frac N\sigma (p-1)}}d\tau ds + \int_0^t \frac 1{g(s)} \left[\int_0^t\left(e^{\Gamma(\tau)}\right)'(\varepsilon \mathcal{H}(\tau)-\mathcal{I}(\tau))d\tau+\mathcal{A}(0)\right]ds.
\end{aligned}
\end{equation}

To estimate the left hand side of the above equation, we introduce the following Lemma, which is proved in \cite{II}.

\begin{lemma}\label{bulemma2}
Let 
\begin{equation*}
\chi_\beta :=\left\{
\begin{aligned}
1, &\quad \mbox{if} \quad \beta = -1,\\
0, &\quad \mbox{if} \quad \beta\in(-1,1).
\end{aligned}
\right.
\end{equation*}
Then, 
\begin{equation}
\begin{aligned}
&\frac{d}{ds}\left\{[\log(s+1)+1]^{\chi_\beta}(\Gamma(s)+1)\int_0^s\frac {g(\tau)e^{\Gamma(\tau)}}{G(\tau)+1}\mathcal{R}(\tau)d\tau\right\}\\
&\geq C\frac 1{g(s)} e^{\Gamma(s)}\mathcal{R}(s)+C \frac 1{g(s)}\int_0^s\frac {g(\tau)e^{\Gamma(\tau)}}{G(\tau)+1}\mathcal{R}(\tau)d\tau ds
\end{aligned}
\end{equation}
\end{lemma}

Thus, applying Lemma \ref{bulemma2}, 
\begin{equation}\label{bu37}
    \begin{aligned}
        &\int_0^t \frac 1{g(s)} e^{\Gamma(s)}\mathcal{R}(s)ds+C \int_0^t \frac 1{g(s)}\int_0^s\frac {g(\tau)e^{\Gamma(\tau)}}{G(\tau)+1}\mathcal{R}(\tau)d\tau ds\\
        \leq & C\int_0^t \left[\frac 1{g(s)} e^{\Gamma(s)}\mathcal{R}(s)+\frac 1{g(s)}\int_0^s\frac {g(\tau)e^{\Gamma(\tau)}}{G(\tau)+1}\mathcal{R}(\tau)d\tau\right] ds\\
        \leq &\int_0^t \frac{d}{ds}\left\{[\log(s+1)+1]^{\chi_\beta}(\Gamma(s)+1)\int_0^s\frac {g(\tau)e^{\Gamma(\tau)}}{G(\tau)+1}\mathcal{R}(\tau)d\tau\right\}ds\\
        =& C [\log(t+1)+1]^{\chi_\beta}(\Gamma(t)+1)\int_0^t\frac {g(\tau)e^{\Gamma(\tau)}}{G(\tau)+1}\mathcal{R}(\tau)d\tau.
    \end{aligned}
\end{equation}
Set
\begin{equation*}
X(t):=[\log(t+1)+1]^{\chi_\beta}(\Gamma(t)+1)\int_0^t\frac {g(\tau)e^{\Gamma(\tau)}}{G(\tau)+1}\mathcal{R}(\tau)d\tau,
\end{equation*}
and
\begin{equation*}
Y(t):=\int_0^t \frac 1{g(s)}\int_0^s \frac {(e^{\Gamma(s)}-e^{\Gamma(\tau)})g(\tau)\mathcal{R}(\tau)^p}{(G(\tau)+1)^{\frac N\sigma (p-1)}}d\tau ds.
\end{equation*}
Combining (\ref{bu35}) and (\ref{bu37}), we obtain
\begin{equation}\label{bu38}
    C_2 X(t)\geq C_3 Y(t)+\int_0^t \frac 1{g(s)} \left[\int_0^t\left(e^{\Gamma(\tau)}\right)'(\varepsilon \mathcal{H}(\tau)-\mathcal{I}(\tau))d\tau+\mathcal{A}(0)\right]ds.
\end{equation}
Now, we estimate the second term of the right hand side in (\ref{bu38}). Recall that
\begin{equation}
    \begin{aligned}
        \mathcal{H}(t)&=\varepsilon^{-1}(\mathcal{C}(t)-\mathcal{E}_2(t))\\
        &=\left(\frac{G(t)+1}{2G(t)+1}\right)^{N/\sigma}\int_{\R^N}\varphi(\frac x{(2G(t)+1)^{1/\sigma}}) (a_0+g(0)a_1)(x)dx\\
        &\ +\left(\frac{G(t)+1}{2G(t)+1}\right)^{N/\sigma} \frac{N g(0)^2}{\sigma (2G(t)+1)}\int_{\R^N}\varphi(\frac x{(2G(t)+1)^{1/\sigma}}) a_0(x)dx\\
        &\ -\left(\frac{G(t)+1}{2G(t)+1}\right)^{N/\sigma}\frac {g(0)^2}{2\sigma (G(t)+1)}\int_{\R^N}\nabla \varphi(\frac x{(2G(t)+1)^{1/\sigma}})\cdot \frac x{(2G(t)+1)^{1/\sigma}} a_0(x)dx
    \end{aligned}
\end{equation}
Note that \(G(t)=\int_0^t g(s)ds\rightarrow \infty\) as \(t\rightarrow +\infty\), since \(C^{-1}b(t)^{-1} \le g(t) \le Cb(t)^{-1}\). Therefore, 
\begin{equation}
    \lim_{t\rightarrow +\infty} \mathcal{H}(t) = 2^{-N/\sigma}\varphi(0)\int_{\R^N}(a_0+g(0)a_1)(x)dx>0,
\end{equation}
by the assumption on the initial data. Then, by an argument similar to that in \cite{II}, we get 
\begin{equation}
    \int_0^t \frac 1{g(s)} \left[\int_0^t\left(e^{\Gamma(\tau)}\right)'(\varepsilon \mathcal{H}(\tau)-\mathcal{I}(\tau))d\tau+\mathcal{A}(0)\right]ds\geq C_3\varepsilon e^{\Gamma(t)},
\end{equation}
for \(t>t_\varepsilon\), which implies
\begin{equation}\label{bu42}
    C_2 X(t)\geq C_3 Y(t)+C_4\varepsilon e^{\Gamma(t)}\quad \mbox{for}\quad t>t_\varepsilon.
\end{equation}
Here, \(t_\varepsilon\) depends polynomially on \(\varepsilon\).

On the other hand, direct calculation implies
\begin{equation}
    g(t)Y''(t)+g'(t)Y'(t)-Y'(t)=\frac 1{g(t)}\int_0^t \frac {e^{\Gamma(\tau)})g(\tau)\mathcal{R}(\tau)^p}{(G(\tau)+1)^{\frac N\sigma (p-1)}}d\tau.
\end{equation}
And by H\"{o}lder's inequality, 
\begin{equation}
    \begin{aligned}
        X(t)&=[\log(t+1)+1]^{\chi_\beta}(\Gamma(t)+1)\int_0^t\frac {g(\tau)e^{\Gamma(\tau)}}{G(\tau)+1}\mathcal{R}(\tau)d\tau\\
        &\leq [\log(t+1)+1]^{\chi_\beta}(\Gamma(t)+1)\left(\int_0^t\frac {g(\tau)e^{\Gamma(\tau)}}{(G(\tau)+1)^{\frac N\sigma(p-1)}}\mathcal{R}(\tau)^pd\tau\right)^{\frac 1p} \left(\int_0^t\frac {g(\tau)e^{\Gamma(\tau)}}{(G(\tau)+1)^{p'-\frac N\sigma}}d\tau\right)^{\frac1{p'}}.
    \end{aligned}
\end{equation}
Therefore, 
\begin{equation}\label{bu45}
    \begin{aligned}
         &g(t)Y''(t)+g'(t)Y'(t)-Y'(t)\\
         &=\frac 1{g(t)}\int_0^t \frac {e^{\Gamma(\tau)})g(\tau)\mathcal{R}(\tau)^p}{(G(\tau)+1)^{\frac N\sigma (p-1)}}d\tau\\
         &\geq \frac 1{g(t)}X(t)^p\left\{[\log(t+1)+1]^{\chi_\beta}(\Gamma(t)+1)\right\}^{-p}\left(\int_0^t\frac {g(\tau)e^{\Gamma(\tau)}}{(G(\tau)+1)^{p'-\frac N\sigma}}d\tau\right)^{-\frac p{p'}}.
    \end{aligned}
\end{equation}
Thus, combining the above inequality with (\ref{bu42}), we get
\begin{equation}\label{bu46}
    \begin{aligned}
         &g(t)Y''(t)+g'(t)Y'(t)-Y'(t)\\
         &\geq \frac 1{g(t)}X(t)^p\left\{[\log(t+1)+1]^{\chi_\beta}(\Gamma(t)+1)\right\}^{-p}\left(\int_0^t\frac {g(\tau)e^{\Gamma(\tau)}}{(G(\tau)+1)^{p'-\frac N\sigma}}d\tau\right)^{-\frac p{p'}}\\
         &\geq C \frac 1{g(t)}(Y(t)+\varepsilon e^{\Gamma(t)})^p\left\{[\log(t+1)+1]^{\chi_\beta}(\Gamma(t)+1)\right\}^{-p}\left(\int_0^t\frac {g(\tau)e^{\Gamma(\tau)}}{(G(\tau)+1)^{p'-\frac N\sigma}}d\tau\right)^{-\frac p{p'}}.
    \end{aligned}
\end{equation}

Set \(W(t):=e^{-\Gamma(t)}Y(t)+\varepsilon\). Then, by direct calculations, we have
\begin{equation}\label{bu47}
    \begin{aligned}
         &W''(t)+b(t)W'(t)\\
         &\geq C \frac {W(t)}{g(t)^2}e^{(p-1)\Gamma(t)}\left\{[\log(t+1)+1]^{\chi_\beta}(\Gamma(t)+1)\right\}^{-p}\left(\int_0^t\frac {g(\tau)e^{\Gamma(\tau)}}{(G(\tau)+1)^{p'-\frac N\sigma}}d\tau\right)^{-\frac p{p'}}.
    \end{aligned}
\end{equation}
We treat the right hand side of the above inequality by the similar method in \cite{II} to get that for \(t>t_\varepsilon\),
\begin{equation}
    \begin{cases}
        (t+1)^{\beta}W''(t)+CW'(t)\geq C\frac{W(t)^p}{t+1},\quad &\mbox{if}\quad \beta\in (-1,1),\\
        (t+1)^{-1}W''(t)+CW'(t)\geq C\frac{W(t)^p}{(t+1)(\log (t+1)+1},\quad &\mbox{if}\quad \beta=-1.
    \end{cases}
\end{equation}

Thus, applying Lemma A.1 and Lemma A.2 in \cite{II} respectively, we can conclude that
\begin{equation}
    T_0\leq 
    \begin{cases}
        \exp{(C\varepsilon^{-(p-1)}},\quad &\mbox{if}\quad \beta\in (-1,1),\\
        \exp{(\exp{(C\varepsilon^{-(p-1)}})},\quad &\mbox{if}\quad \beta=-1.
    \end{cases}
\end{equation}

\section*{Acknowledgements}
The first author is supported by Science Foundation of Zhejiang Sci-Tech University (ZSTU) under Grant No.13062120-Y. The second author is supported by JST CREST Grant Number JPMJCR1913, Japan and Grant-in-Aid for Young Scientists Research (No.19K14581), Japan Society for the Promotion of Science.


\begin{thebibliography}{99}

\bibitem{AD} H.S. Aslan and T.A. Dao, On the Cauchy problem for semilinear $\sigma$-evolution equations with time-dependent damping, Math. Meth. Appl. Sci. (2023), 1–38, DOI 10.1002/mma.9857.

\bibitem{DL} M. D'Abbicco and S. Lucente, A modified test function method for damped wave equations. Adv. Nonlinear Stud. 13 (2013), no. 4, 867-892.

\bibitem{DLR} M. D'Abbicco, S. Lucente and M. Reissig, Semi-linear wave equations with effective damping. Chin. Ann. Math. Ser. B 34 (2013), no. 3, 345-380.

\bibitem{DP}  T.A. Dao and N.H. Son, Critical curve for a weakly coupled system of semi-linear $\sigma$-evolution equations with frictional damping, Rocky Mountain Journal of Mathematics, 52 (2022), 299-321.

\bibitem{DR} T.A. Dao and M. Reissig, Blow-Up Results for Semi-Linear Structurally Damped $\sigma$-Evolution Equations. In: Cicognani, M., Del Santo, D., Parmeggiani, A., Reissig, M. (eds) Anomalies in Partial Differential Equations. Springer INdAM Series, vol 43 (2021), Springer, Cham, 213-245, http://doi.org/10.1007/978-3-030-61346-4$\_$10.

\bibitem{FIW2019} K. Fujiwara, M. Ikeda and Y. Wakasugi, Estimates of lifespan and blow-up rates for the wave equation with a time-dependent damping and a power-type nonlinearity. Funkcial. Ekvac. 62 (2019), no. 2, 157-189.

\bibitem{FIW2020} K. Fujiwara, M. Ikeda and Y. Wakasugi, The Cauchy problem of the semilinear second order evolution equation with fractional Laplacian and damping, Nonlinear Differ. Equ. Appl. 28 (2021), no. 6, Paper No. 63, https://doi.org/10.1007/s00030-021-00723-6.

\bibitem{FO1}K. Fujiwara and T. Ozawa, Lifespan of srong solutions to the periodic nonlinear Schr\"{o}dinger equation without gauge invariance, J. Evol. Equ., 17 (2017), 1023-1030.

\bibitem{FO2}K. Fujiwara and T. Ozawa, Finite timeblowup of solutions to the nonlinear Schr\"{o}dinger equation without gauge invariance, J. Math. Phys., 57 (2016), 082103.

\bibitem{FRRSY} I. Friedman, O.Ria\~{n}o, S. Roudenko, D. Son and K. Yang, Well-posedness and dynamics of solutions to the generalized KdV with low power nonlinearity, Nonlinearity, 36 (2023), 584-635.

\bibitem{Hao} C. Hao, Lecture notes on introduction to Harmonic analysis, Chinese Academy of Sciences, 2012.

\bibitem{II} M. Ikeda and T. Inui, The sharp estimate of the lifespan for the semi-linear wave equation with time-dependent damping, Differential Integral Equations 32 (2019)
1-36.

\bibitem{IO} M. Ikeda, T. Ogawa, Lifespan of solutions to the damped wave equation with a critical nonlinearity, J. Differential Equations 261(3) (2016) 1880-1903.

\bibitem{ISW} M. Ikeda,  M. Sobajima and Y. Wakasugi, Sharp lifespan estimates of blowup solutions to semi-linear wave equations with time-dependent effective damping. J. Hyperbolic Differ. Equ. 16 (2019), no. 3, 495-517.

\bibitem{IW} M. Ikeda and Y. Wakasugi, A note on the lifespan of solutions to the semilinear damped wave equation, Proc. Amer. Math. Soc. 143 (2015), no. 1, 163-171.

\bibitem{KL} K. Kim and S. Lim, Asymptotic behaviors of fundamental solution and its derivatives to fractional diffusion-wave equations, J. Korean Math. Soc. 53 (2016), No. 4, 929-967.

\bibitem{LaiZhou} N. A. Lai and Y. Zhou, The sharp lifespan estimate for semi-linear damped wave
equation with Fujita critical power in higher dimensions, J. Math. Pures Appl. 123 (2019) 229-243

\bibitem{LZ} T. Li and Y. Zhou, Breakdown of solutions to $\square u+u_t=u^{1+\alpha}$, Discrete Cont. Dynam. Syst., 1 (1995), 503-520.

\bibitem{N2003}  K. Nishihara, $L^p-L^q$ estimates of solutions to the damped wave equation in 3-dimensional space and their application, Math. Z. 244 (2003) 631-649.

\bibitem{LI} J. Lin and M. Ikeda, Weighted Fractional Commutator Estimate and Application to the Fractional Wave Equations with Time-dependent Damping, preprint.

\bibitem{NLZ} J. Lin, K. Nishihara and J. Zhai, Critical exponent for the semi-linear wave equation with time-dependent damping, Discrete Contin. Dyn. Syst. 32 (2012) 4307-4320.

\bibitem{N}K. Nishihara, Asymptotic behavior of solutions to the semilinear wave equation with time-dependent damping. Tokyo J. Math. 34 (2011), no. 2, 327-343.

\bibitem{QD} Y. Qian, T. Dao, On the Cauchy problem for a weakly coupled system of semi-linear $\sigma$-evolution equations with double dissipation, arXiv: 2311.06663v1.

\bibitem{TY1} G. Todorova and B. Yordanov, Critical exponent for a nonlinear wave equation with damping, C. R. Acad. Sci. Paris I, 330 (2000), 557-562.

\bibitem{TY2} G. Todorova and B. Yordanov, Critical exponent for a nonlinear wave equation with damping, J. Differential Equations, 174 (2001), 464-489.

\bibitem{W} Y. Wakasugi, Scaling variables and asymptotic profiles for the semilinear damped wave equation with variable coefficients, J. Math. Anal. Appl., 447 (2017), pp. 452-487.

\bibitem{Z} Qi S. Zhang, A blow-up result for a nonlinear wave equation with damping:the critical case, C.R. Acad. Sci. Paris I, 333 (2001),109-114.

\end{thebibliography}
\end{document}